\input amstex
\input psfig.tex

\input psfig

\newif\ifboxfigure      
\boxfigurefalse

\def\BoxIt#1#2{
	\vbox{\hrule
	\hbox{\vrule\kern#2\vbox{\kern#2#1\kern#2}\kern#2\vrule}
		   \hrule}}

\def\insertRaster #1 pixels #2  by #3 scaled #4 {
			\medskip
			 \hbox to \hsize{%

			 \hss
			 \RasterBox {#1} {#2} {#3} {#4}
			 \hss
			 }%
}

\def\RasterBox #1 #2 #3 #4{


\dimen5=65pt
\divide\dimen5 by 72

\dimen0=#2\dimen5
\divide\dimen0 by 1000
\dimen1=#3\dimen5
\divide\dimen1 by 1000
\dimen2=#3\dimen5
\divide\dimen2 by 1000
\dimen3=#2\dimen5
\divide\dimen3 by 1000

\setbox4=\hbox to #4\dimen0{
 \vbox to #4\dimen1{
 \vss
 \psfig{figure=#1,height=#4\dimen2,width=#4\dimen3}
 }
 \hss
 }
 \ifboxfigure\BoxIt{\box4}{0pt}
 \else\box4
 \fi
 }


\documentstyle{amsppt}


\define \C{{\Bbb C}}
\define \CC{\overline \C} 
\define \R{{\Bbb R}}
\define \Z{{\Bbb Z}}
\define \CZ{{\Bbb C}/{\Bbb Z}}

\redefine \t{\tau}

\define \N{\Cal N}

\define \s{\sigma}
\define \Q{\Cal Q} 
\redefine \Re{\operatorname{Re} }
\redefine \Im{\operatorname{Im} }
\define \Dom{{\Cal Dom} }
\define \cE{\Cal E}
\define \tE{\tilde {\Cal E}}
\define \hE{\hat {\Cal E}}
\define \cR{\Cal R}
\define \tR{\tilde {\Cal R}}
\define \hR{\hat {\Cal R}}

\define\bM{\partial M}
\define\HD{\operatorname{H-dim} }
\define\hd{\operatorname{hyp-dim} }

\redefine\l{\lambda}
\redefine\L{\Lambda}
\redefine\epsilon{\varepsilon}
\predefine\vphi{\phi}
\redefine\phi{\varphi}
\redefine\emptyset{\vphi}

\topmatter
\title
The Hausdorff dimension of the boundary of the Mandelbrot set and 
Julia sets 
\endtitle

\author
Mitsuhiro Shishikura
\endauthor

\rightheadtext\nofrills{Hausdorff dimension of the boundary of 
the Mandelbrot set and 
Julia sets} 

\affil 
Tokyo Institute of Technology \\
and \\
State University of New York at Stony Brook \\ 
 \\
\endaffil

\comment 
\address
Tokyo Institute of Technology \\ 
Department of Mathematics \\
Ohokayama, Meguro, Tokyo 152, Japan \\ 
email mitsu@math.titech.ac.jp \\ 
and \\ 
Institute for Mathematical Sciences \\
State University of New York at Stony Brook \\ 
Stony Brook, NY 11794-3660, USA
\endaddress
\email mitsu@math.sunysb.edu \endemail
\endcomment

\dedicatory
Dedicated to Prof. John W. Milnor on the occasion of his sixtieth birthday
\enddedicatory

\comment 

\keywords 
Mandelbrot set, Julia set, Hausdorff dimension, 
parabolic periodic point, Ecalle cylinder
\endkeywords

\subjclass
\endsubjclass

\endcomment

\abstract 
It is shown that the boundary of the Mandelbrot set $M$ has Hausdorff dimension two 
and that for a generic $c \in \bM$, the Julia set of $z \mapsto z^2+c$ also has 
Hausdorff dimension two.  The proof is based on the study of the bifurcation 
of parabolic periodic points.  
\endabstract 

\endtopmatter

\NoRunningHeads 


\document 

\heading 
Introduction
\endheading 

\medskip 

The dynamics of complex quadratic polynomials $P_c(z)=z^2+c$ 
has been studied extensively recent years (for example [DH]).  
The main interest in this subject is to study the nature of the Julia sets 
$J_c$ in the dynamical plane and the Mandelbrot set $M$ in the parameter 
space.  The boundary of $M$ 
also has a meaning as ``the locus of bifurcation", or more precisely, by 
Ma\~n\'e-Sad-Sullivan [MSS] or by Lyubich [Ly1], 
$\bM=\{c \in \C|\ P_c \text{ is not J-stable }\} $  
(see \S 1 for definition).   
In this paper, we are concerned with 
the Hausdorff dimension (denoted by $\HD(\cdot)$) of these sets.  
Some of the consequences are:  

\medskip 

\proclaim{Theorem A}
$$\HD(\bM)=2.$$
Moreover for any open set  $U$  which intersects 
$\bM$,  we have  $\HD(\bM \cap U)=2$. 
\endproclaim

\medskip 

\proclaim{Theorem B} For a generic $c \in \bM$, 
$$\HD J_c=2. $$   
In other words, there exists a residual (hence dense) subset $\Cal R$ of $\bM$ 
such that if $c \in \Cal R$,  then 
$\HD J_c=2$.  
\endproclaim 
\medskip 

\proclaim{Theorem C}
There exists a residual set $\Cal R'$ of $\R/\Z$ such that if $P_c$ has a 
periodic 
point with multiplier $\exp(2\pi i \alpha)$  with $\alpha \in \Cal R'$, then 
$\HD J_c =2.$  
\endproclaim

 Theorem A was conjectured by many people (for example, Mandelbrot [Ma], 
Milnor [Mi2, Conjecture 1.7]).  
It means that the bifurcation locus is large in dimension, and this explains 
the complexity of $M$, which is observed by many numerical experiments.  
As to the Julia sets, 
if $P_c$ is hyperbolic, or if 
$0$ is strictly preperiodic, $\HD J_c$ is less than 2 (see \S 1, Property (1.4)).
However it was conjectured that there exists a sequence of parameters 
such that $\HD J_c$ tends to 2.  Theorem B gives a stronger solution  
to this conjecture.  
We will see that the method in this paper applies 
to other families under certain condition.  

\medskip 

The above theorems are obtained as consequences of Theorems 1 and 2 stated 
in \S 1.  Theorem 1 amounts to comparing the Hausdorff dimension of the set of 
J-unstable parameters with that of a certain subset (``hyperbolic subset'') 
of the Julia set.  
It reflects the similarity between the Mandelbrot set and some Julia sets.  
(Compare with Tan Lei [T].)  
Theorem 2 is the most important result in this paper, and it 
assures that one can obtain maps whose Julia sets have Hausdorff dimension 
(or ``hyperbolic dimension'') arbitrarily close to 2, from ``the 
secondary bifurcation'' of a parabolic periodic point.  
See Figures 1, 2 and the remark after Theorem 2 in \S 2.  

\midinsert 
\hbox to \hsize{ 
\hfil 
\hbox{\RasterBox {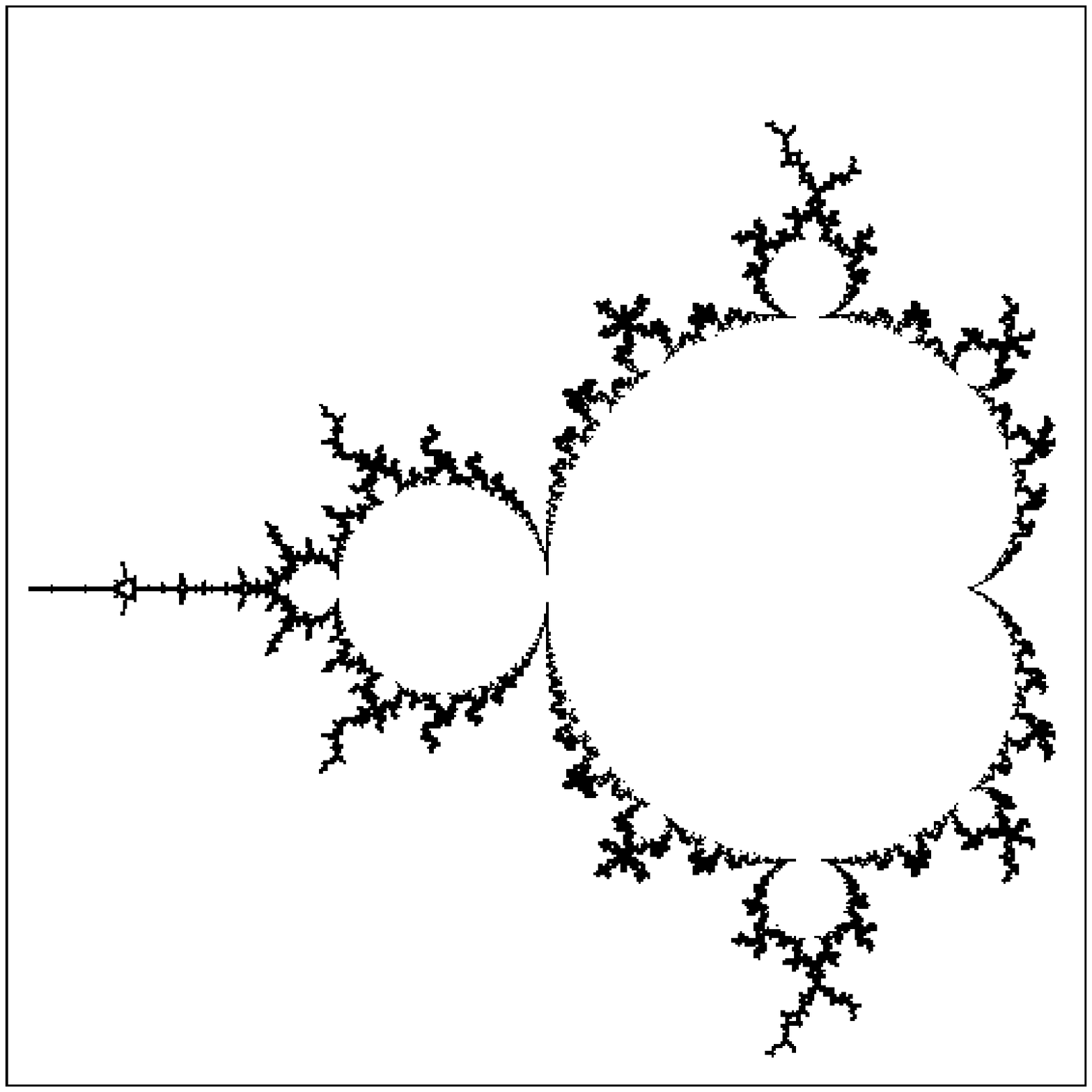} {480} {480} {270} }
\hfil 
\hbox{\RasterBox {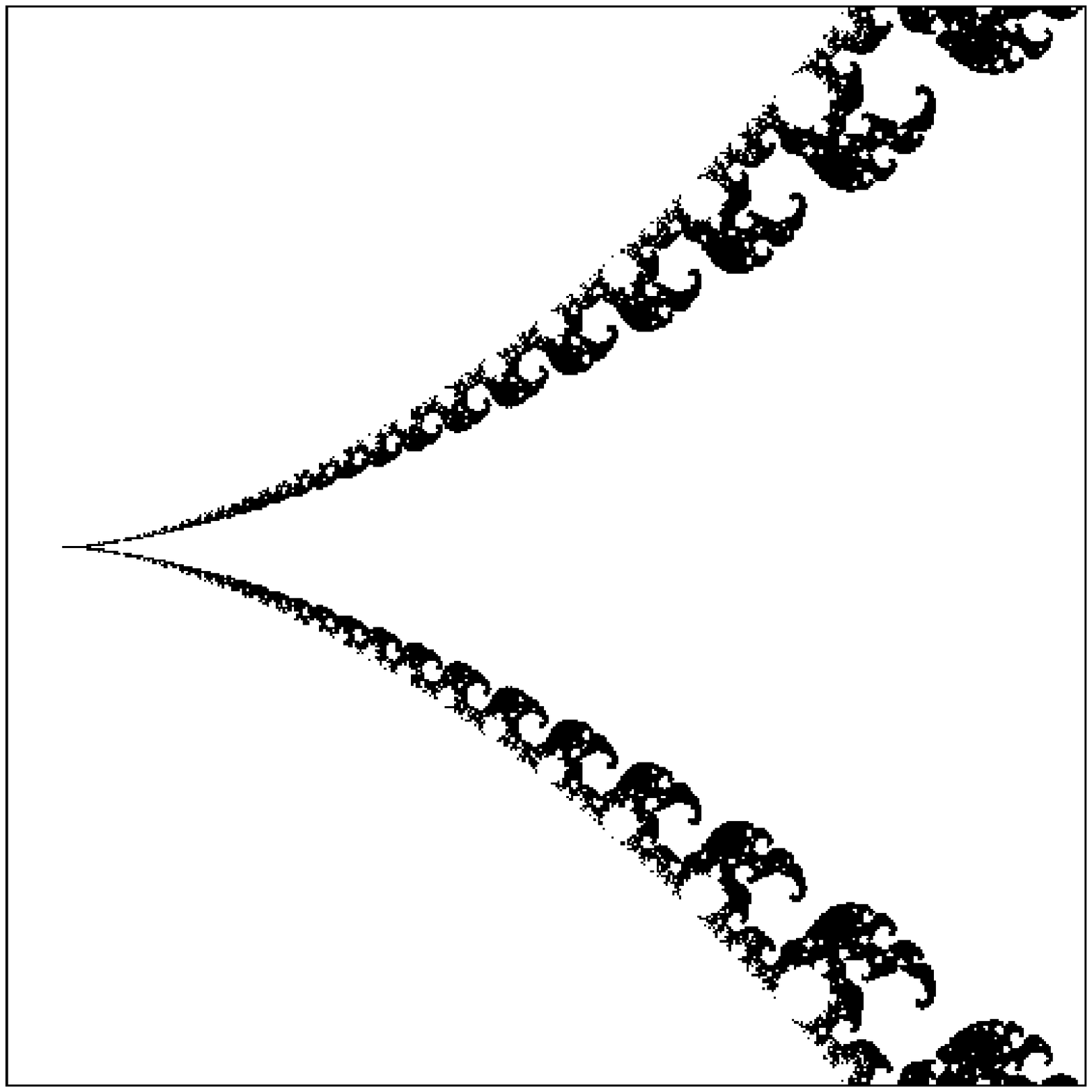} {480} {480} {270} } 
\hfil 
\hbox{\RasterBox {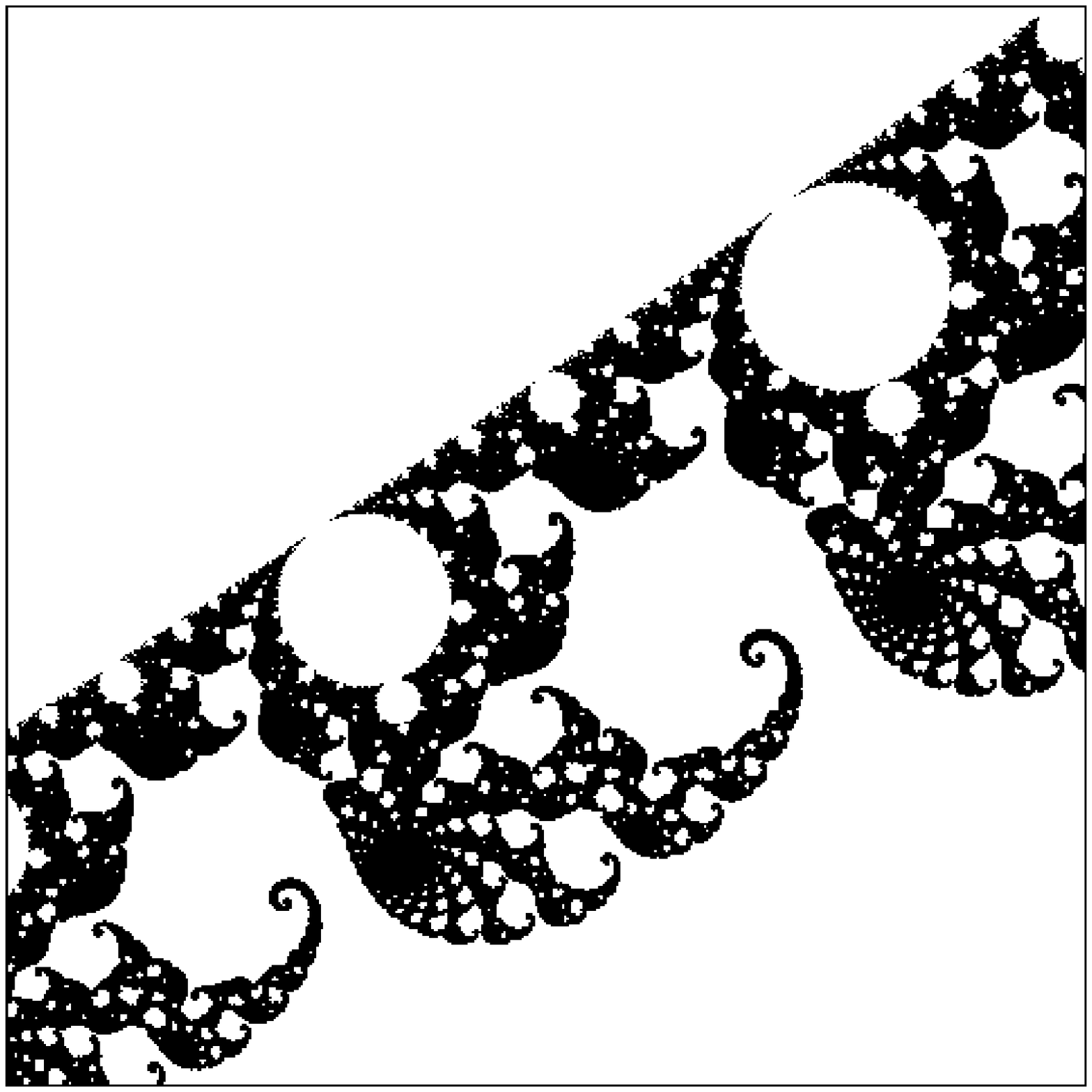} {480} {480} {270} }
\hfil 
}
\hbox to \hsize{ 
\hfil 
\hbox{\RasterBox {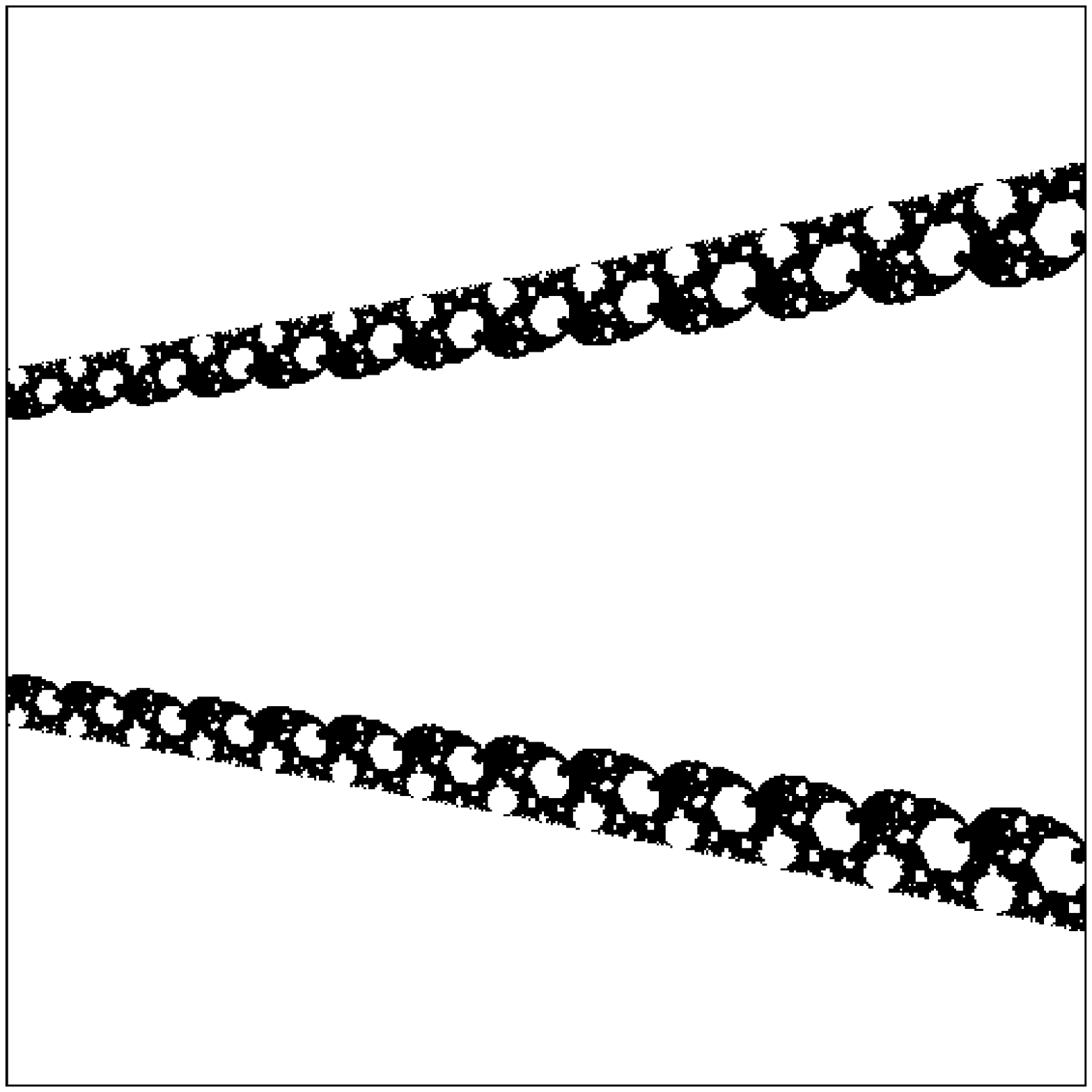} {480} {480} {270} }
\hfil 
\hbox{\RasterBox {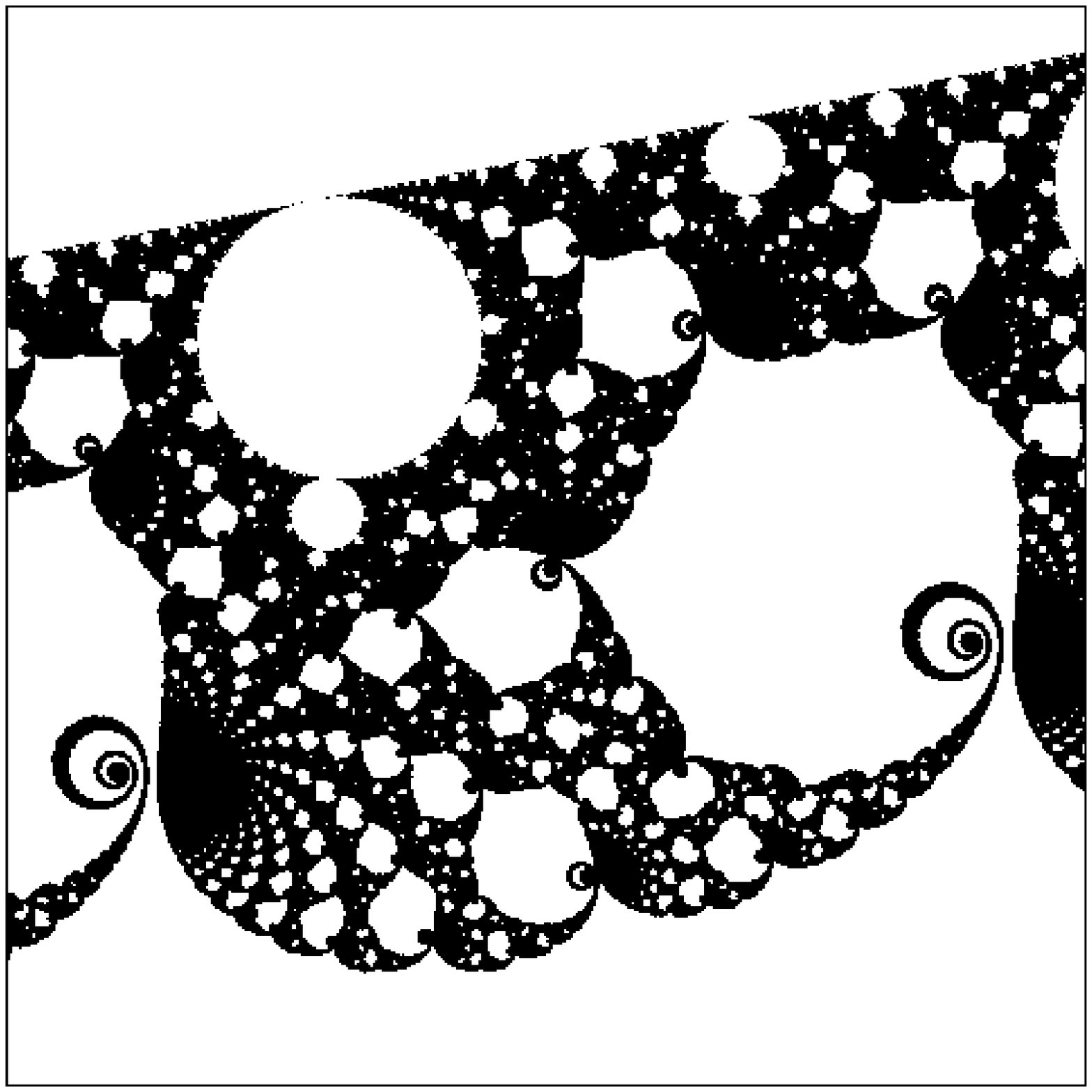} {480} {480} {270} } 
\hfil 
\hbox{\RasterBox {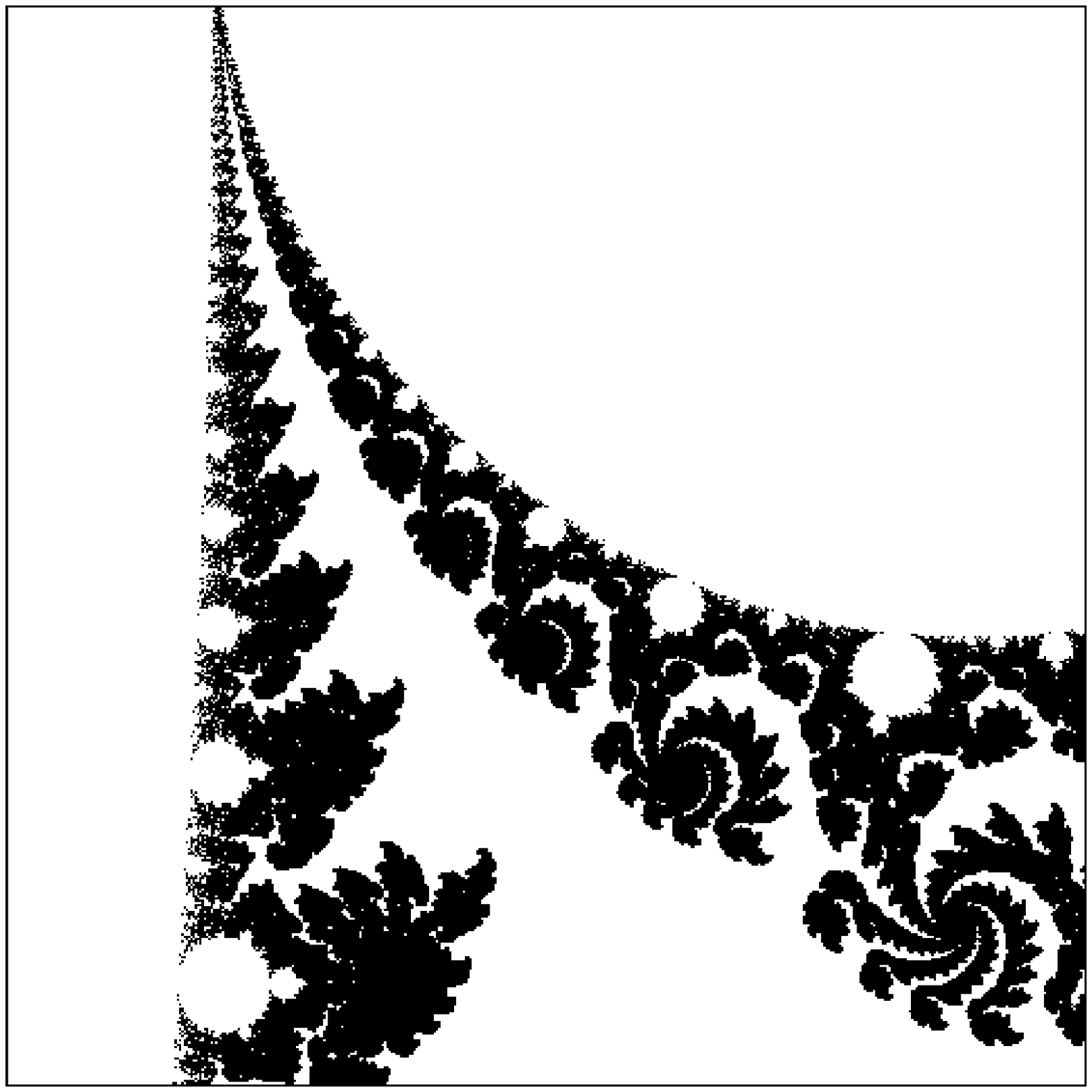} {480} {480} {270} }
\hfil 
}
\botcaption{Figure 1} 
The boundary of he Mandelbrot set (top left) and its 
blow-ups near the cusp $c=\frac 14$.  
\endcaption 
\endinsert 

\midinsert 
\hbox to \hsize{ 
\hfil 
\hbox{\RasterBox {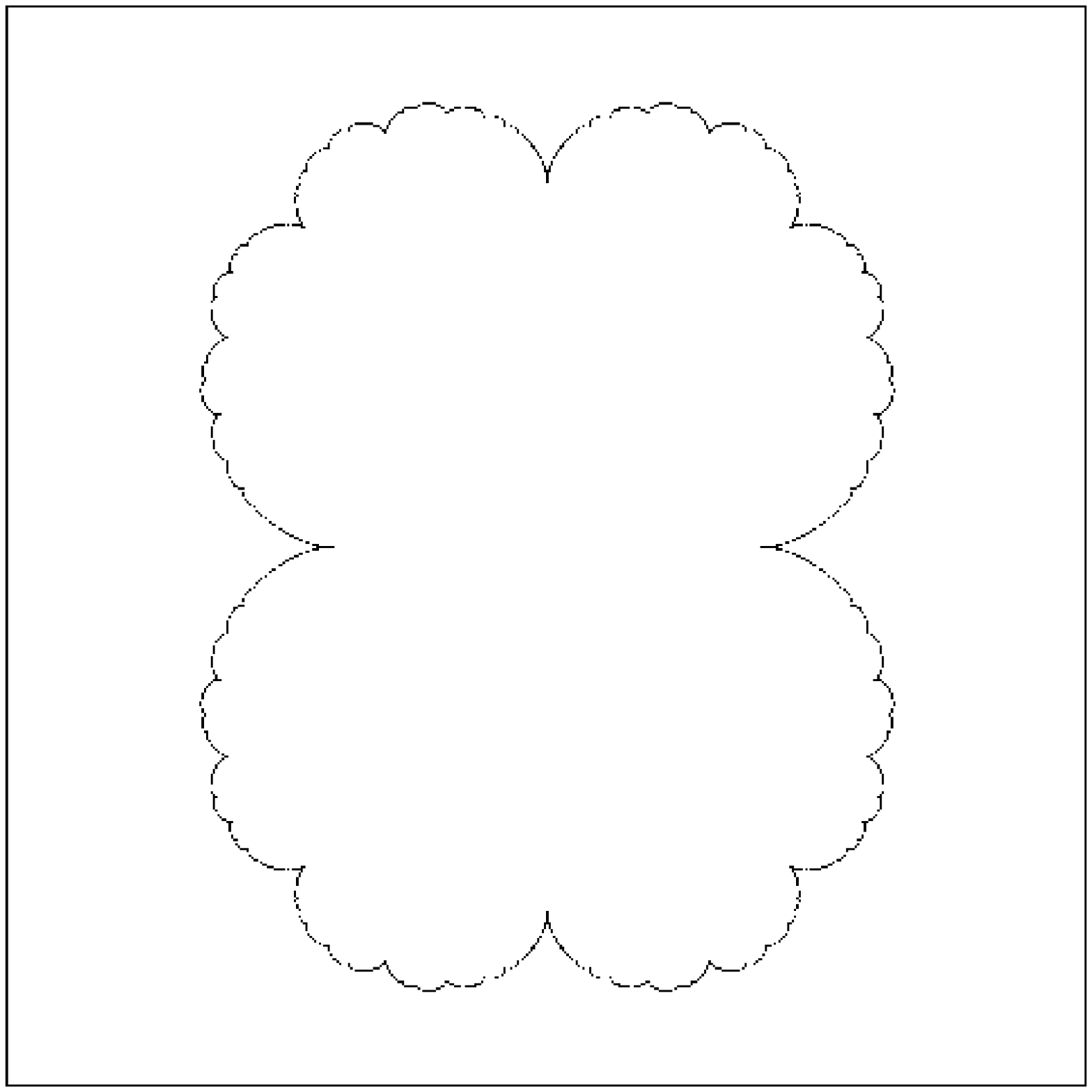} {480} {480} {270} }
\hfil 
\hbox{\RasterBox {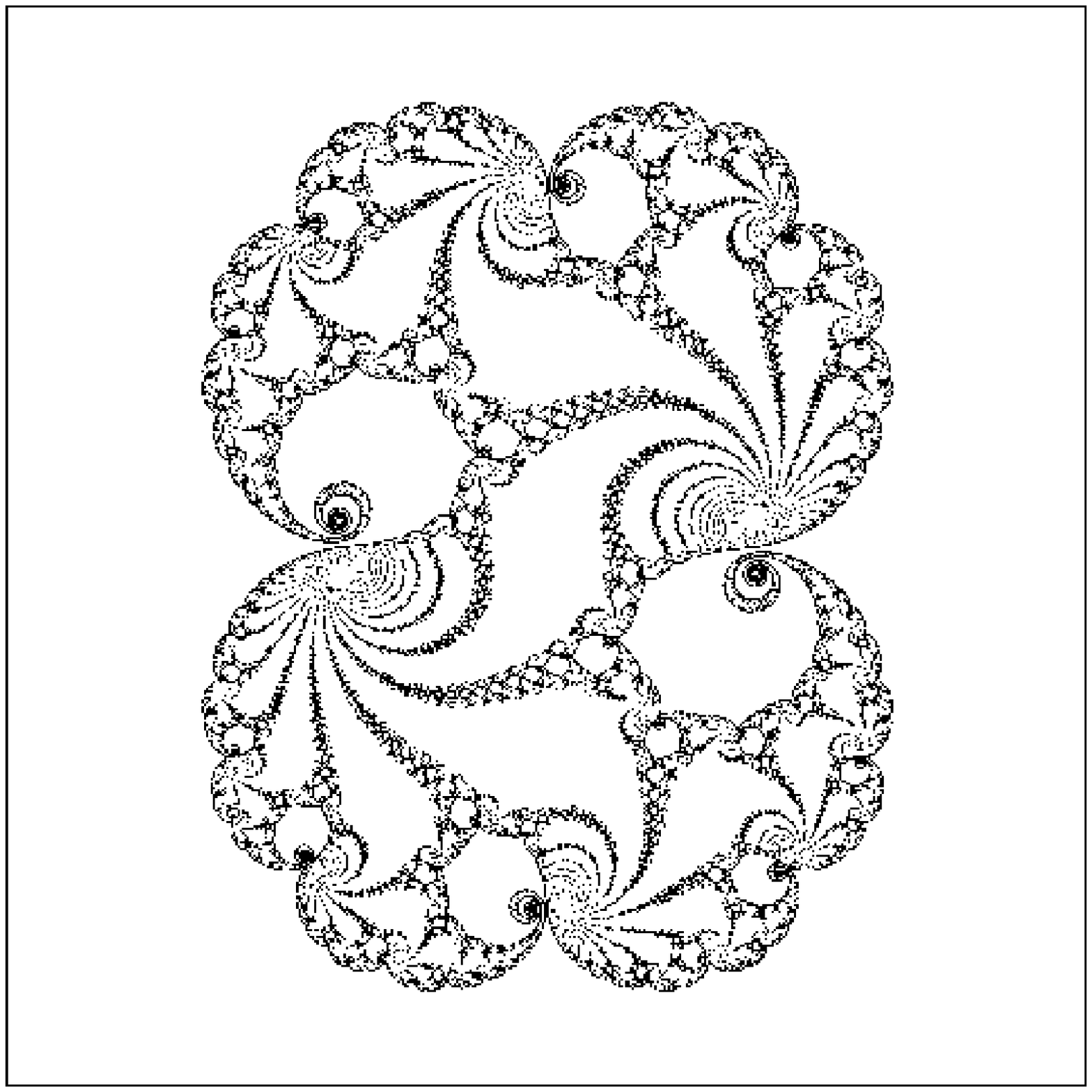} {480} {480} {270} } 
\hfil 
\hbox{\RasterBox {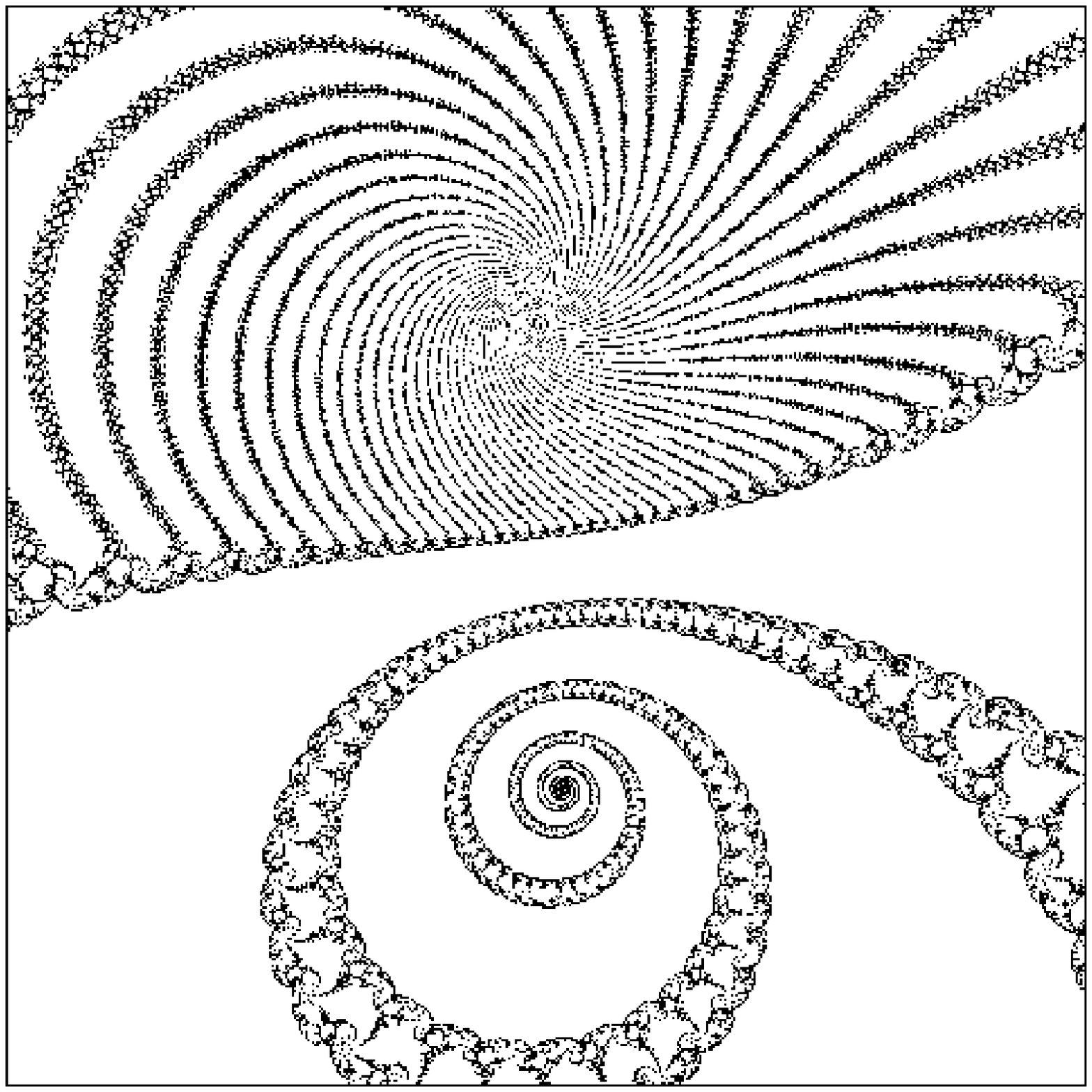} {480} {480} {270} }
\hfil 
}
\botcaption{Figure 2}
The Julia sets for $c=\frac 14 $ (left) and for $c=0.25393+0.00048 i$ (middle)  
and its blow up near the fixed points (right)
\endcaption 
\endinsert

It was already observed that after a small perturbation of a parabolic 
periodic point, the Julia set may inflate drastically.  In fact, the proof of 
Theorem 2 shows that such an inflated part of the Julia set can have 
Hausdorff dimension close to 2.   
The main tool in the study of such a bifurcation is the theory of Ecalle 
cylinders, which was introduced by Douady-Hubbard [DH] and developed by 
Lavaurs [L].  Using the Ecalle cylinder, we introduce  
a new renormalization procedure associated with parabolic fixed points 
(see Remark (ii) in \S 6).  
Our result can be interpreted as follows: The renormalization induces 
a map between old and new dynamical planes, which resembles an exponential 
map.  Comparing this observation with McMullen's result [Mc] which claims that 
the Julia set of an exponential map has always Hausdorff dimension 2, 
it becomes conceivable that a certain subset of the Julia set 
can have Hausdorff dimension close to 2.  
The proof in this paper will justify this argument, or even more, 
twice renormalization is enough in order to attain dimension two.   

\medskip

  One can compare the above theorems with Jakobson's result for the 
family of unimodal interval maps [J], M. Rees' result for a certain family 
of rational maps [Re] and Benedicks-Carleson's 
result on the family of H\'enon maps [BC].  
These results show the existence of a ``chaotic dynamics'' 
for a set of positive Lebesgue measure of parameters.  
For example, M. Rees' result shows that there exists a set of positive 
Lebesgue measure of parameters for which the Julia sets are the whole 
Riemann sphere.  
Such parameters are found near a special parameter for which all critical 
points are strictly preperiodic.  
For this parameter,  
the map has good ergodic theoretical properties and the Julia set 
is the whole sphere.  

On the other hand, for a polynomial $P_c$ acting on $\C$, the Julia set $J_c$ 
or the filled-in Julia set $K_c$ can never be the entire plane, since 
there is always the basin of $\infty$.  So there always exists some 
orbits which escape to $\infty$ from the neighborhood of $K_c$.  
Moreover, as remarked above, if the critical point $0$ is strictly preperiodic, 
$\HD J_c <2$.  For example if $c=-2$, $J_{-2}=[-2,2] \subset \R$ and 
$\HD J_{-2}=1$.  
Therefore one can hardly expect an analogous result 
or approach for the family $P_c$ as those 
of Jakobson and Rees.  Instead, in this paper, 
we use the perturbation of parabolic periodic points.  

\medskip 

As for the area (the 2-dimensional Lebesgue measure), it is conjectured that 
$\bM$ and $J_c$ (for any $c$) have area zero.  There are partial results: 
the set of parameters in $\bM$ for which $P_c$ are not infinitely renormalizable 
has area zero [Sh].  If $P_c$ has no irrationally indifferent periodic point 
and is not infinitely renormalizable, then the Julia set $J_c$ has area zero 
[Ly2] and [Sh].  

\medskip 

This paper is organized as follows: 
In Section 1, we define the notion of hyperbolic subsets and the hyperbolic 
dimension, state our main results (Theorems 1 and 2, Corollary 3) and 
assuming these results, we prove 
Theorems A, B and C.  
We prove basic properties of hyperbolic subsets and hyperbolic dimension 
in Section 2.  Theorem 1 is proved in Section 3, using holomorphic motions.  
The rest of the paper is devoted to the proof of Theorem 2.  The theory of 
the parabolic 
bifurcation and Ecalle cylinders is reviewed in Section 4.  Further 
properties of the Ecalle transformation are studied in Section 5.  
After these preparations, Theorem 2 is proved in Section 6 (the case with 
multiplier 1) and in Section 7 (the other cases).  In the Appendix, we give 
the proof of the facts stated in Sections 4 and 7.  

\bigskip 

\definition{Acknowledgement} 
I would like to thank specially Curt McMullen for introducing me these problems and 
having inspiring discussions with me throughout on this subject; 
A. Douady for his lectures which introduced me the theory of Ecalle cylinders; 
and also A. Hinkkanen, M. Lyubich, M. Rees, D. Sullivan and J. Milnor 
for helpful discussions and comments while this paper was written.  
Finally I am grateful to the Institute for Mathematical Sciences, State University 
of New York at Stony Brook for the hospitality.  

Computer pictures have been produced using J. Milnor's program.  
\enddefinition  
\newpage

\comment 
\midinsert 
\centerline{\hbox{
\psfig{figure=,height=5cm,width=5cm} 
\psfig{figure=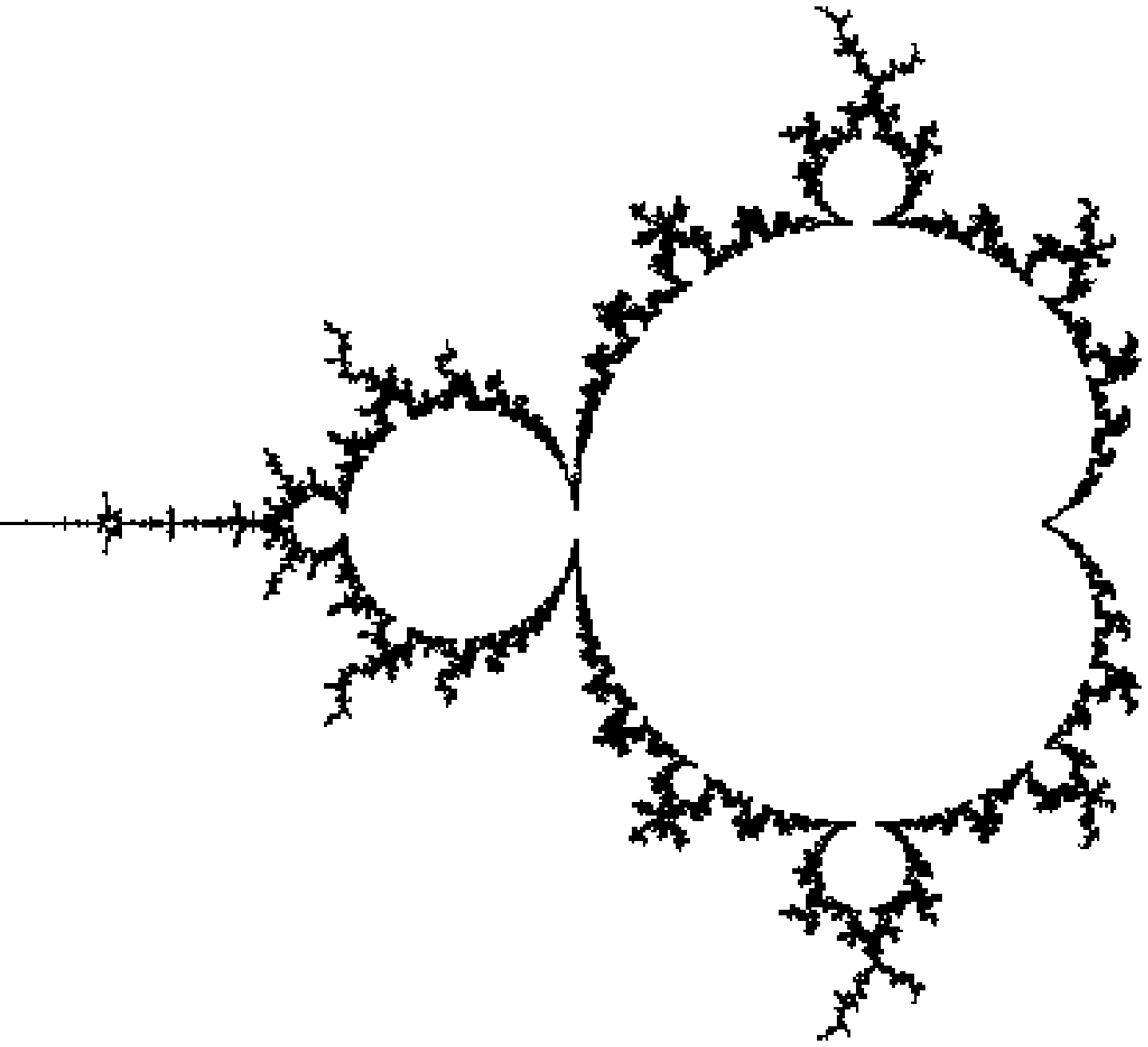,height=5cm,width=5cm} 
}} 
\botcaption{Figure 1 and 2}
the Mandelbrot set 
\endcaption 
\endinsert 

  One can compare the above theorems with Jakobson's result for the 
family of unimodal interval maps [J] and M. Rees' result for a certain family 
of rational maps [Re].  These results  state that 
in the parameter space, there exists a set of positive 
Lebesgue measure, consisting of parameters for which 
the maps exhibit a ``chaotic dynamics'' on the 
whole phase space (the interval or the Riemann sphere), 
for example, in the case of rational maps, the Julia set is the whole sphere.  
They find the set of parameters as above 
near a special parameter for which all critical 
points (only one in the unimodal map case) are strictly preperiodic.  
Hence for this parameter,  
the map has good ergodic theoretical properties.   Remark that 
Benedicks-Carleson's 
result on the family of H\'enon maps [BC] also relies on the perturbation 
of such a unimodal map into two dimensional maps.

\endcomment

\heading
\S 1. Some definitions and Main results.  
\endheading

\medskip 

  In this section, we define the hyperbolic subset and the hyperbolic dimension 
for a rational map, and state our main results-- Theorems 1 and 2.  
Assuming these results, we give the proof of Theorems A, B and C in Introduction.  
\definition{Definitions}  Let  $f$  be a rational map.  A closed subset 
$X$ of $\CC$  is 
called a {\it hyperbolic subset} for $f$, if 
$f(X) \subset X$ and there exist positive constants $c$ and $\kappa>1$ 
such that 
$$\|(f^n)'\| \geq c \kappa^n \text{ on } X \ \text{  for } n \geq 0, $$
where  $\|\cdot \|$ denotes 
the norm of the derivative with respect to the spherical metric of $\CC$.  
(A similar notion was already seen in [Ru] and [Ly1].  
The condition requires that $f$ is expanding on $X$, 
however we call such a set ``hyperbolic'' following the standard terminology 
for dynamical systems.)  

The {\it hyperbolic dimension} of $f$ is
$$\hd(f)=\sup\{\HD(X)|\ X \text{ is a hyperbolic subset for } f \}.$$

\enddefinition

\proclaim{
Properties of hyperbolic subsets and hyperbolic dimension} 
\newline 
Let $X$ be a hyperbolic subset of $f$. 
\newline 
{\rm \bf (1.0)} There are no critical points of $f$ in $X$.  
\newline 
{\rm \bf (1.1)} $X$ is a subset of the Julia set $J(f)$ of $f$.     
Hence 
$$\HD J(f) \geq \hd(f).$$   
{\rm \bf (1.2)} $X$ is ``stable" under a perturbation, i.e. there exists a 
neighborhood 
$\Cal N$ of $f$ in the space of rational maps of the same degree, 
such that if $g \in 
\Cal N$ then $g$ has a hyperbolic subset $X_g$ and there is a homeomorphism 
$\iota_g:X \to X_g$ which conjugates $f$ to $g$.  
Moreover, for each $z \in X$,  
$\iota_g(z)$ is a complex analytic function in $g$, and $\iota_f=id_X$.  
($\{\iota_g\}$ is a holomorphic motion in the sense of \S 3.)  
\newline  
{\rm \bf (1.3)} $f \mapsto \hd(f)$ is lower semi-continuous, or equivalently, 
for any number c, the set $\{f|\ \hd(f)>c \}$ is open.  
\newline 
{\rm \bf (1.4)} Suppose that $f$ is a hyperbolic rational map (i.e., 
all critical points are attracted to 
attracting periodic orbits) or a subhyperbolic polynomial (every critical point  
is either attracted to an attracting periodic orbit or pre-periodic, see [DH]).  
Then 
$$\HD J(f)=\hd(f).$$  
Moreover  
$J(f)$ has positive and finite $\delta$-dimensional Hausdorff 
measure if $\delta=\HD J(f)$,  and $\HD J(f)<2$.  
\endproclaim

For the proof and remarks, see \S 2.  

\remark{Problem} 
When does $\hd(f)$ coincide with $\HD J(f)$?
\endremark

\medskip

\definition{Definition}
A family $\{f_\l|\ \l \in \L\ \}$ of rational maps is {\it J-stable at} 
$\l_0 \in \L$, if there exists a continuous map $h: \L ' \times J(f_{\l_0}) 
\to \CC$, such that 
$\L'$ is a neighborhood of $\l_0$ in $\L$, $h_\l \equiv h(\l,\cdot)$ 
is a conjugacy from $(J(f_{\l_0}),f_{\l_0})$ 
to $(J(f_\l),f_\l)$  and  $h_{\l_0}=id_{J(f_{\l_0})}$.  We also say that 
$f_{\l_0}$ is J-stable in this family, if there is no confusion.  
\enddefinition

\medskip

\proclaim{Theorem 1}
  Let $\{f_\l|\ \l \in \L\ \}$ be a complex analytic 
family of rational maps of degree $d\ (>1)$, where $\L $ is an open set in 
$\C$.  Suppose $f_{\l_0}$ $(\l_0 \in \L)$ is not J-stable in this family.  
Then 
$$
 \HD \biggl\{  \l \in \L \bigg| \matrix 
f_\l \ \text{\rm is not J-stable and has a hyperbolic subset} \\
 \text{\rm containing a forward orbit of a critical point } 
\endmatrix \biggr\} 
\geq \  \hd(f_{\l_0}).  $$

\endproclaim
The proof will be given in \S 3.  
\bigskip

\definition{Definition}
A periodic point is called {\it parabolic} if its multiplier is a root of 
unity.  The {\it parabolic basin} of a parabolic periodic point $\zeta$ 
of period $k$ is 
$$\{z \in \CC | \ f^{nk} \to \zeta \ (n \to \infty) \ 
\text{ in a neighborhood of } z \ \}, $$ 
and the {\it immediate parabolic basin} of $\zeta$ is 
the union of periodic connected components of the parabolic basin.
\enddefinition

\proclaim{Theorem 2}
Suppose that a rational map $f_0$ of degree $d\ (>1)$ has a parabolic fixed 
point $\zeta$ with multiplier 
$\exp(2\pi i p/q)$ $(p, q \in \Z,\ (p,q)=1)$ and that 
the immediate parabolic basin of $\zeta$ 
contains only one critical point of $f_0$.  Then: 
\newline 
for any $\epsilon > 0$ and $b>0$, there exist a neighborhood $\Cal N$ of 
$f_0$ in the space of rational maps of degree $d$, a neighborhood 
$V$ of $\zeta$ in $\CC$, positive integers $N_1$ and $N_2$ such that 
if $f \in \Cal N$, and if $f$ has a 
fixed point in $V$ with multiplier  $\exp(2\pi i \alpha)$, where  
$$
q \alpha = p  \pm \cfrac 1 \\
\         a_1 \pm \cfrac 1 \\
\         a_2 + \beta \ \endcfrac 
$$
with integers $a_1 \geq N_1$, $a_2 \geq N_2$ and 
$\beta \in \C, \ 0\leq \Re \beta <1, 
|\Im \beta | \leq b$, \ then 
$$\hd(f) > 2-\epsilon.$$

\endproclaim 

The proof will be given in \S\S 6 and 7, after preparations in \S\S 4 and 5.  

\bigskip

Figure 3 shows the region for $\alpha$ described in Theorem 2.  

\midspace{7cm} \caption{Figure 3. the region for $\alpha$}    

The condition $a_1 \geq N_1$ is in fact unnecessary, since $a_1$ must be large 
when $f$ is close to $f_0$.  However if we take a family $\{f_\alpha \}$ such that 
$f_\alpha(0)=0$, $f'_\alpha(0)=\exp(2 \pi i \alpha)$, then the condition 
in Theorem 2 can be expressed in terms of $\alpha$ (without $\N$ and $V$).  

The condition on $\alpha$ has a meaning that $f$ is obtained 
by a ``secondary bifurcation'', in the following sense:  
The primary bifurcation of $f_0$ 
produces a sequence of maps having parabolic fixed points
with multipliers $\exp(2\pi i \alpha_n)$, where 
$q\alpha=p \pm \frac 1n \ (n =1,2, \dots )$.  Then the bifurcation from 
these maps gives rise to the maps as in 
Theorem 2.  

\bigskip
  There are immediate consequences of Theorems 1 and 2 for the family 
$P_c$.  Before stating them, let us recall the definition of $J_c$ and $M$:   
$$\aligned 
K_c &= \{z \in \C | \ P_c^n(z) \not \to \infty \text{  as } n \to \infty \} 
\ \ (\text{the {\it filled-in Julia set} of } P_c), \\
J_c &= \partial K_c = \text{the closure of } \{  
\text{ repelling periodic points of } P_c\} 
\ (\text{the {\it Julia set} of } P_c), \\ 
M &=\{ c \in \C| \ 0 \in K_c \} = \{c \in \C| \ K_c \text{ is connected }\} 
\ \ (\text{the {\it Mandelbrot set}}). 
\endaligned $$

\proclaim{Corollary 3}
(i) If $U$ is an open set containing $c \in \bM$, then 
$$\HD(\bM \cap U) \geq \HD\{ c \in \bM \cap U| \ 0 \text{ is non-recurrent under } 
P_c \} \geq \hd(P_c).$$
(ii)If $P_c$ has a parabolic periodic point with multiplier 1, then 
there exists a sequence $\{c_n\}$ in $\bM$  such that $c_n \to c$ and 
$\hd(P_{c_n}) \to 2, \ \text{as}\ n\to \infty$.  
\endproclaim

\demo{Proof of Corollary 3}  (i) Immediate from Theorem 1 and the fact that 
$c \in \bM$ if and only if $P_c$ is J-unstable ([MSS],[Ly1]).  
\newline 
(ii) It is known that that if $P_c$ has a parabolic periodic point 
$\zeta$ of order $k$, 
then $f_0=P_c^k$ satisfies the hypothesis of Theorem 2, since there is 
only one critical point (in $\C$) for $P_c$.  
Any parabolic periodic point of $P_c$ is not persistent (otherwise all $P_c$ 
would have parabolic periodic points), and it can be 
perturbed into a periodic point whose multiplier is as in 
Theorem 2.  If $\Im \beta=0$, then the new periodic point is indifferent, 
and the perturbed polynomial is also J-unstable.  
One can thus obtain the sequence $\{c_n\}\ (\subset \bM)$.  Furthermore, it 
is also possible to choose $c_n$ such that $0$ is strictly preperiodic 
under $P_{c_n}$, since such parameters are 
dense in $\bM$ by [MSS].  
\qed \enddemo 

Now we can prove the Theorems A, B and C.  

\demo{\bf Proof of Theorem A}  By Ma\~n\'e-Sad-Sullivan [MSS] or Lyubich [Ly1], 
parameters for which $P_c$ has (non-persistent) parabolic periodic points 
are dense in $\bM$.  The assertion follows immediately from Corollary 3.  
\qed \enddemo 
   
\demo{\bf Proof of Theorem B} 
Let $\Cal R_n=\{c\in \bM|\ \hd(P_c)>2-\frac 1n \}$ $(n=1,2\cdots)$.  
Then $\Cal R=\bigcap_{n\geq0} \Cal R_n = \{c \in \bM |\ \hd(P_c)=2\ \} \subset 
\{c \in \bM |\ \HD\ J(P_c)=2\ \}$.  By Property (1.3), $\Cal R_n$ are open in 
$\bM$.  Moreover $\Cal R_n$ is dense in $\bM$, by Corollary 3(ii) and the above 
remark (the density of parabolic parameters in $\bM$).  
Hence $\Cal R$ is residual.   
(A residual set is a set containing the intersection of a countable 
collection of open dense subsets of $\bM$.)  Such an $\Cal R$ is dense 
in $\bM$  by Baire's theorem.  
\qed \enddemo 

\demo{\bf Proof of Theorem C} 
Let $W$ be a hyperbolic component, i.e. a connected component of the set of 
$c$'s such that $P_c$ has an attracting periodic point.  Then it was shown by 
Douady-Hubbard [DH] that there exists a homeomorphism from $\overline W$ onto 
the closed unit disc (conformal in $W$) defined by the multiplier of a non-
repelling periodic orbit.  So in $\partial W$, the  parameters with 
parabolic periodic points are dense.  
By the proof of Corollary 3(ii) and a similar argument 
to the proof of Theorem B, we can prove that for generic $\alpha \in \R/\Z$, 
if $c \in \partial W$ and $P_c$ has a periodic 
point with multiplier $\exp(2\pi i \alpha)$, then 
$\HD J_c =2.$  However, there are only countably many hyperbolic components.  
Hence the assertion follows.  
\qed \enddemo 

\medskip 

\remark{Remark}  
(i) It also follows easily from Theorem 2 that 
$$\sup_{c \in W_0} \HD(P_c) = \sup_{c \in \C \setminus M} \HD(P_c) = 2\ , $$
where $W_0=\{c|\ P_c \text{ has an attracting fixed point }\}$ ({\it the  
cardioid}).  
\newline 
(ii)  It is easy to see that a similar result holds for other families 
of rational maps which have 
``only one critical point'' that can be involved in the parabolic basin.  
For example, $f_a(z)=z^3+az^2$, $g_b(z)=(z^2+b)/(z^2-1)$.   
\endremark

\bigskip 
\bigskip

\heading
\S 2. Hyperbolic subsets and hyperbolic dimension 
\endheading

\medskip 

In this section, we give the proof of the Properties (1.0)-(1.3).  
We also give an example of the hyperbolic subset and an estimate of 
its Hausdorff dimension, which will be used later.  

\demo{\bf Proof of the properties of hyperbolic subsets and 
hyperbolic dimension}  
\newline 
(1.0) is obvious.  
\newline 
(1.1) The family $\{f^n\}$ cannot be normal in any open set 
intersecting $X$, since the derivatives grow exponentially.  Hence the 
assertion follows.  
\newline 
(1.2) This is a well-known fact.   The outline of the proof is as follows.  
There exists a neighborhood $V$ of $X$ such that for $g$ near $f$, 
$g$ is expanding on $V$.  Hence for $x \in X$,  the orbit $\{f^n(x)\}_{n=0}^\infty$ 
with respect to $f$ is a pseudo-orbit for $g$, which can be traced by 
an orbit (a real orbit) of a point $y$ for $g$.  
(The Pseudo-Orbit Tracing Property, see [Bo1].)  Then let $y= \iota_g(x)$.  
In fact, $y$ can be expressed as 
$$
y=\lim_{n \to \infty} g_{x_0}^{-1} \circ g_{x_1}^{-1} \circ \cdots 
\circ g_{x_n}^{-1}(x_{n+1}), 
$$
where $x_j=f^j(x)$ and $g_z^{-1}$ is the branch of $g^{-1}$ defined in a 
neighborhood of $f(z)$ such that $g_z^{-1}(f(z))$ is near $z$.  
By the expanding property for $g$ near $f$,  
$\iota_g$ is well-defined and conjugates 
$f|_X$ to $g|_{X_g}$ with $X_g=\iota_g(X)$.  Moreover, $\iota_g(x)$ depends 
analytically in $g$, since the convergence in the above is uniform.  

(In [Ly1], the analyticity is proved under the assumption that periodic points are 
dense in $X$, which is in fact unnecessary by the above.)  

\comment 
Let $X$ be a hyperbolic subset of $f$ and $m$ a positive integer such that 
$||(f^m)'||>2$  on $X$.  There exist a neighborhood $\N$ of $f$, a neighborhood 
$V$ of $X$ in $\CC$  and $\delta>0$  such that if $g \in \N$ and if $x,y \in V$ 
with $d(x,y)<\delta$, then 
$$d(g^m(x),g^m(y))\geq 2 d(x,y). \tag{*}$$
Then there exists an $\epsilon>0$  
($\epsilon < \delta/2$) such that if $x \in X$ and $y \in \CC$ satisfy 
$d(f^m(x),y)<\epsilon$, then there exists a unique $z \in V$ such that 
$d(x,z)<\delta/2$ and $g^m(z)=y$.  Denote this $z$ by $\theta_g(x,y)$.  
It can be easily seen that $\theta_g(x,y)$ depends continuously on $y$, 
analytically on $g$ and locally constant with respect to $x$.   
Let us take $\N$ smaller so that $d(f^m(x),g^m(x))<\epsilon$.  It follows 
from (*) that 
if $d(f^m(x), y)<\epsilon$, then 
$$d(x, \theta_g(x,y)) \leq \frac 12 d(g^m(x),y) \leq 
\frac 12 (d(g^m(x),f^m(x))+d(f^m(x),y)) <\epsilon.$$     
Hence we can define inductively $\theta_g^k(x,y)$ for $y$ 
satisfying $d(f^{mk}(x), y)<\epsilon$ ($k\geq 1$), by 
$$\theta_g^k(x,y)=\theta_g(x,\theta_g^{k-1}(f^m(x),y)).  $$
Let $\iota_g^{(k)}(x)=\theta_g^k(x,f^{mk}(x))$ ($x \in X$), then 
$\iota_g^{(k)}(x)=\theta_g(x,\iota_g^{(k-1)}(f^m(x)))  $ 
$g^m \iota_g^{(k)}(x)= \iota_g^{(k-1)} f^m(x)$  
and 
$$
\aligned 
d(\iota_g^{(k+1)}(x),\iota_g^{(k)}(x)) &\leq 
\frac 12 d(g^m \iota_g^{(k+1)}(x) ,g^m \iota_g^{(k)}(x)) = 
\frac 12 d(\iota_g^{(k)}(x), \iota_g^{(k-1)}(x)) \\
&\leq \cdots 
\leq \left(\frac 12 \right)^k d(\iota_g^{(1)}(x), x) 
\leq \left(\frac 12 \right)^k \epsilon .
\endaligned $$ 
Hence there exists the limit 
$$\iota_g(x)=\lim_{k \to \infty} \iota_g^{(k)}(x)$$ 
uniformly in $x\in X$ and $g\in \N$.  
Since $\iota_g^{(k)}$ varies continuously in $x$ and analytically in $g$, 
so does $\iota_g(x)$.  Obviously $g^m \iota_g= \iota_g f^m$, $\iota_f(x)=x$ 
and $d(x, \iota_g(x))\leq \epsilon$.  If $\iota_g(x)=\iota_g(y)$, 
then $\iota_g(f^{mk}x)=\iota_g(f^{mk}y)$ hence $d(f^{mk}x, f^{mk}y)\leq 
2 \epsilon <\delta$ ($k\geq 1$).  By (*), we have $x=y$, hence 
$\iota_g$ is injective.  
Moreover if we take $\N$ smaller so that  
$d(\iota_g(fx),g\iota_g(x)) <\delta$ and $\iota_g(X) \subset V$, then 
the distance between $g^{mk}\iota_g(fx)=\iota_g f(f^{mk}(x)) \in V$ and 
$g^{mk}(g\iota_g(x))= g\iota_g(f^{mk}(x)) \in V$ is less than $\delta$ 
($k \geq 1$).    
It follows from (*) that $\iota_g f(x)=g\iota_g(x)$.   
\qed 
\endcomment 

\noindent 
(1.3) Let $\Cal N$, $X_g$ and $\iota_g$ be as in (1.2).  (Suppose 
$\Cal N$ is open.)  It is enough to 
prove that  
$\Cal N \ni g \mapsto \HD X_g$ is continuous.  
We will prove this fact in \S 3, using a result on holomorphic motions.  
It is also possible to prove it directly.  
In fact, one can estimate the exponent of the H\"older continuity 
of $\iota_g$, in the proof of (1.2).  
\newline 
(1.4) We do not use this fact 
for the proof of our main theorems.  If $f$ is hyperbolic, then 
$J(f)$ itself is a hyperbolic subset, hence 
$\HD J(f)=\hd(f)$; 
the second assertion follows from Bowen's formula ([Bo], [Ru]); and 
the fact $\HD J(f)<2$ can be shown by a standard argument using 
the expanding property of $f$ 
and the Lebesgue's density theorem. (See [Su].)   
The case with preperiodic critical points will be discussed in another paper, 
but the fact $\HD J(f) <2$ seems to be already known.

\bigskip 

In the proof of Theorem 2 (\S 6), we will construct a special kind of 
hyperbolic sets which are described as follows.  

  Suppose that $U$ is a simply connected open set of $\CC$ 
(with $\sharp(\CC-U) \geq 2$); $U_1, \dots , U_N $ 
are disjoint simply connected open subsets of $U$ with 
$\overline U_i \subset U$; $n_1, \dots, n_N$ are positive integers 
such that $f^{n_i}$ maps $U_i$ onto  $U$ bijectively  
$(i=1,\dots,N)$.  It follows from Schwarz' lemma that $\tau_i \equiv 
(f^{n_i}|_{U_i})^{-1} : U \to U_i$ is a contraction with respect to the 
Poincar\'e metric of $U$ (at least on $\cup_i \overline U_i$).  
So there exists a Cantor set $X_0$ generated by $\{\tau_i\}$, that is,   
$X_0$ is the minimal non-empty closed set satisfying 
$$X_0=\tau_1(X_0) \cup \dots \cup \tau_N(X_0).$$

\medskip

\proclaim{Lemma 2.1} 
The set  
$X=X_0 \cup f(X_0) \cup \dots \cup f^{M-1}(X_0), $ 
(where  $M=\max n_i$) is a hyperbolic subset of $f$.   
\endproclaim

\demo{Proof}  Obviously $X$ is closed. 
We have $f(X) \subset X$, since 
$f(f^{M-1}(X_0))=f^M(\tau_1(X_0) \cup \cdots \cup \tau_N(X_0))
=f^{M-n_1}(X_0) \cup \cdots \cup f^{M-n_N}(X_0) \subset X .$

Note that $f^{n_i}|_{U_i}$ is expanding on $X_0 \cap U_i$ with respect to 
the Poincar\'e metric of $U$ which is equivalent to the spherical metric on 
$X_0$.  
In order to see that $f$ is expanding on $X$, it suffices to factorize  
the iterate $f^n$ at $z \ (\in X)$  to 
$$
f^{j_1}\circ (f^{n_{i_1}}|_{U_{i_1}}) \circ \dotsc 
\circ (f^{n_{i_k}}|_{U_{i_k}}) \circ f^{j_2}
$$
where $i_1,\dots, i_k \in \{1,\dots,N\}$ and $0\leq j_1, j_2 <M$ are determined 
by  $z'=f^{j_2}(z) \in X_0$, 
$z' \in \tau_{i_k} \circ \cdots \circ \tau_{i_1}(X_0)$ and 
$n=j_1+n_{i_1}+\dotsc+n_{i_k}+j_2$.  
\qed (Lemma 2.1) \enddemo 

\bigskip
We only use the simplest estimate for the Hausdorff dimension of $X_0$.  
Suppose $\infty \notin U$.

\proclaim{Lemma 2.2}
Let $\delta=\HD X_0$.  Then 
$$
1 \geq \sum_{i=1}^N \inf_U |\tau'_i|^\delta \geq 
N \cdot \left(\max_i\  \sup_{U_i} |(f^{n_i})'|\right)^{-\delta},  
$$
hence
$$
\delta \geq \frac{\log N}{\log \left({\underset i \to \max}\ 
{\underset {U_i} \to \sup} |(f^{n_i})'| \right)} \ .
$$
\endproclaim

\medskip

\demo{Proof}  The first inequality is well-known; it can be proven, for example, 
from Bowen's formula ([Bo2] and [Ru]).  The rest is immediate.  
\qed (Lemma 2.2) 
\enddemo

\bigskip 
\bigskip

\enddemo
\heading
\S 3. Holomorphic motions. 
\endheading

\medskip

\definition{Definition}  Let $X$ be a subset of $\CC$ and $\L$ a complex 
manifold with a base point $\l_0 \in \L$.  A family 
of maps $i_\l:X \to \CC \ (\l \in \L)$  is called 
a {\it holomorphic motion} of $X$, if each $i_\l$ is injective,  
$i_{\l_0}=identity_X$ and for each $z \in X$, $i_\l(z)$ is analytic in 
$\l$.  We also say that $X_\l \equiv i_\l(X)$ is a holomorphic 
motion of $X$.  We are mostly interested in the case 
$\L=\{ \l \in \C\bigl| \ |\l|<R\} \ (R>0)$  
with the base point $\l_0=0$.  
\enddefinition

\proclaim {Lemma 3.1}
If \ $i_\l:X \to \CC $ $(\ |\l| < R)$ is a holomorphic motion, then both 
$i_\l$ and $i_\l^{-1}$ are H\"older continuous with exponent $\alpha(|\l|/R)$, 
where  $\alpha: (0,1) \to \R_+$ is a function (independent of the motion) 
satisfying 
$\alpha(t) \nearrow 1$, as $t \searrow 0$.  
\endproclaim

\demo{Proof}  The improved $\l$-lemma in [ST] (see also [MSS], [BR]) 
implies that $i_\l$ can be extended to a $K(|\l|/R)$-quasiconformal mapping 
and $K(t) \searrow 1 \ (\text{as } t \searrow 0)$.  Since $K$-quasiconformal 
mapping is $1/K$-H\"older continuous (Mori's inequality, see [A]),  
the assertion holds with $\alpha(t)=1/K(t)$.  For example, by [BR], one can have 
$$1 \geq \alpha(t) \geq (1-3t)/(1+3t), \ \text{for}\  0<t<1/3 .$$
\qed (Lemma 3.1) \enddemo 

As an application, we can now complete the proof of (1.3) (see \S 2), 
by showing that: 
\newline 
for $\Cal N$ (an open neighborhood of $f$) and $X_g$ in (1.2), 
$\Cal N \ni g \mapsto \HD X_g $ is continuous.  

\demo{Proof}  
For any $g_0 \in \Cal N$, $\iota_g \circ \iota_{g_0}^{-1}(z)$ is a holomorphic 
motion of $X_{g_0}$.  Lemma 3.1 implies that for $g$ near $g_0$, 
$\iota_g \circ \iota_{g_0}^{-1}$ is $\alpha'$-bi-H\"older 
$(0<\alpha' =\alpha'(g_0,g) \leq 1)$ and $\alpha' \to 1$ when $g \to g_0$.  
Then we have 
$$
\alpha' \cdot \HD X_{g_0} \leq \HD X_g \leq {\alpha'}^{-1} \cdot \HD X_{g_0}.  
$$ 
Hence $\HD X_g$ is continuous in $g$.   
\qed \enddemo 

\medskip 

\proclaim {Lemma 3.2}
Let  $i_\l:X \to \CC \ \ (\l \in \Delta \equiv \{\l \in \C:\bigl|\ |\l|<1\})$ 
be a holomorphic motion.  
Suppose $v:\Delta \to \CC$ is an analytic map such that 
$v(0)=z_0 \in X $ and $v(\l) \not \equiv i_\l(z_0)$.  Then 
$$
\HD\{\l \in \Delta | \ v(\l) \in i_\l(X) \ \} 
\geq \lim_{r \to 0} \HD (X \cap D_r(z_0)), 
$$
where $D_r(z_0)$ denotes the disc of radius $r$ centered at $z_0$ with respect to 
the spherical metric.   
\endproclaim

\demo{Proof}  Changing the coordinate by M\"obius transformations depending 
analytically on $\l$, we may assume that $z_0=0$ and $i_\l(0) 
\equiv 0$.  First suppose that 
$v'(0) \neq  0.$  There exist positive constants $\rho \ (<1)$ and $a$ such that 
in $\{\l\ \bigl|\  |\l|<\rho\} $, 
$v(\l)$ is injective and $a|\l| \leq |v(\l)| < \infty$.  Let 
$$b_r=\sup\{ \ |i_\l(z)|\ \bigl| \ z\in X \cap D_r(0),  \ |\l| \leq \rho \}. 
$$
It follows from $\l$-lemma [MSS] (or Montel's theorem) 
that $b_r \to 0$  as  $r \to 0$.  
So there is $r_0>0$ such that $a\rho > b_r$ for $0<r<r_0$.  Take such an $r$.  

For $z \in X \cap D_r(0)$ and $|\mu| < R_r \equiv a\rho / b_r$, let 
us consider the equation 
$$
v(\l)-i_{\l \mu}(z)=0  \ \ 
(\l \in \Delta_\mu \equiv \{\l \bigl|\  |\l|< \min\{\rho, \rho /|\mu|\}\}). 
\tag{3.3}
$$
Both $v(\l)$ and $i_{\l \mu}(z)$ are analytic in $\Delta_\mu$, and 
for $\l \in \partial \Delta_\mu \ $  we have 
$$|v(\l)| \geq a\cdot 
\min\{\rho, \rho /|\mu|\} > b_r \geq |i_\l(z)|.$$  
Since $v=0$ has the 
unique solution $0$ in $\Delta_\mu$, 
the equation (3.3) also has a unique solution by Rouch\'e's theorem, and 
it depends analytically on $\mu$.  Moreover, for the same $\mu$ and a 
different $z$, the equation gives a different solution, because of the 
injectivity of $i_\l$.  

Now define 
$$
Y_\mu^r = \{\l \in \Delta_\mu |\ 
v(\l)=i_{\l \mu}(z) \text{ for some } z \in X \cap D_r(0) \ \}. 
$$  
Then by the above, $Y_\mu^r\ (|\mu|<R_r)$  is a holomorphic motion of $Y_0^r$, 
and the injections $j_\mu^r:Y_0^r \to Y_\mu^r$ are given by the following: \  
$\l = j_\mu^r(v^{-1}(z))$ is the unique solution of the equation (3.3).  
Note that 
$Y_0^r=v^{-1}(X \cap D_r(0))$,  hence $\HD Y_0^r=\HD (X \cap D_r(0))$ 
and that 
$Y_1^r \subset \{\l \in \Delta | \ v(\l) \in i_\l(X) \ \}$.  
By Lemma 3.1, $j_\mu^r:Y_0^r \to Y_\mu^r$ is $\alpha(|\mu|/R_r)$-bi-H\"older, 
therefore we have  
$$
\HD \{\l | v(\l) \in i_\l(X) \} \geq \HD Y_1^r \geq 
\alpha(1/R_r) \cdot \HD Y_0^r = \alpha(1/R_r) \cdot \HD (X \cap D_r(0)).
$$  
Letting $r \to 0$, we obtain the desired inequality.  

Let us consider the case $v'(0)=0$.  By the assumption, $v \not \equiv 0$.  
By a coordinate change, we may assume that $\infty \in X$ and 
$i_\l(\infty) \equiv \infty$.  Let $m=order(v,0)$  and  $G(z)=z^m$.  
Define $\tilde X_\l=G^{-1}(X_\l)$ and $\tilde X=\tilde X_0$.  
By lifting $v$ and $i_\l$ by  $G$ which is a branched covering 
branched over $0$ and $\infty$, 
one obtains an analytic map $\tilde v: \L \to \CC$  satisfying 
$v=G \circ \tilde v$, $\tilde v'(0) \neq 0$,  and a holomorphic motion 
$\tilde i_\l: \tilde X \to 
\tilde X_\l $ satisfying $i_\l \circ G = G \circ \tilde i_\l$.  
Hence the inequality holds for $\tilde i_\l$ and $\tilde v$.  On the other hand, 
we have 
$$\{\l | v(\l) \in X_\l \} = 
\{\l | \tilde v(\l) \in \tilde X_\l \} 
\text{  and  } 
\HD (X \cap D_r(0))=\HD (\tilde X \cap G^{-1}(D_r(0))), $$
since $G$ is locally Lipschitz except at $0$ and $\infty$.  
Thus we obtain the inequality for $i_\l$ and $v$.  
\qed (Lemma 3.2) 

\bigskip

Now we can give:

\demo{\bf Proof of Theorem 1}   
For any $\epsilon > 0$, there exists 
a hyperbolic subset $X$ for $f_{\l_0}$ such that 
$\HD X > \hd(f_{\l_0})-\epsilon$.  
By the compactness of $X$, there exists a point $z_0 \in X$ such that  
$\  \lim_{r \to 0} \HD (X \cap D_r(z_0)) = \HD X .$

By Property (1.2) of the hyperbolic subset,  
there exist a neighborhood $\L' (\subset \L)$ of $\l_0$ and 
a holomorphic motion $i_\l:X \to X_\l$ for $\l \in \L'$ 
such that $i_\l \circ f_{\l_0} = f_\l \circ i_\l$, 
$X_\l$ is a hyperbolic set of $f_\l$ and  $i_{\l_0}=id_X$.  
Moreover $\L'$ can be chosen smaller so that 
$\  \lim_{r \to 0} \HD (X_\l \cap D_r(i_\l(z_0))) > 
\HD X - \epsilon$ 
for $\l \in \L'$ , by Lemma 3.1, and that 
the critical points of $f_\l$ do not bifurcate in $\L'-\{\l_0\}$.  

It follows from Ma\~n\'e-Sad-Sullivan's theory (Lemma III.2 [MSS]) that 
there exist $\l_1 \in \L'-\{\l_0\}$, an integer $N>0$ and 
a critical point $c$ of $f_{\l_1}$ such that 
$f_{\l_1}^N(c)=i_{\l_1}(z_0)$.  Then there exists a branch of 
critical points $c_\l$ of $f_\l$ in a neighborhood 
$\L'' (\subset \L')$ of $\l_1$ with $c_{\l_1}=c$, 
(hence $c_\l$ is a meromorphic function).  Note that in the above, 
$\l_1$ and $c$ can be chosen so that $f_\l^N(c_\l)\not \equiv i_\l(z_0)$ 
in $\L''$.  Applying Lemma 3.2 to $i_\lambda\ (\l \in \L'' )$ and 
$v(\l)=f_\l^N(c_\l)$, (after a suitable affine change of parameter),  
one obtains
$$\HD\{\l \in \L'' | \ f_\l^N(c_\l) \in X_\l \ \} \geq 
\lim_{r \to 0} \HD (X_{\l_1} \cap D_r(i_{\l_1}(z_0))) > 
\hd(f_{\l_0}) - 2\epsilon .$$
It is easy to see that if $f_\l^N(c_\l) \in X_\l$, $f_\l$ is not J-stable 
in the family, 
since $f_\l^N(c_\l) \not \equiv  i_\l(z)$ for any $z$.  As $\epsilon > 0$ 
was arbitrary, the theorem is proved.  
\qed (Theorem 1) \enddemo 

\medskip 

\remark{Remark}  
If we assume a certain transversality condition for the motion of the critical 
value $v(\l)$ relative to the hyperbolic set, it is also possible to 
prove a similar result (to Theorem 1 or to Lemma 3.2) without assuming the 
analytic dependence of $i_\l$ and $v$ on the parameter.  
\endremark 

\bigskip
\bigskip

\enddemo
\heading
\S 4. Parabolic bifurcation and Ecalle cylinders 
\endheading

\subheading{4.0 Overview} 

  Let us consider a holomorphic mapping 
$$f_0(z)=z+a_2 z^2 +\dots $$
defined near $0$ with $a_2 \neq 0$.  The origin $z=0$ is a parabolic 
fixed point of $f_0$.  If we perturb $f_0$, this fixed point 
bifurcates into two fixed points near $0$ in general.  

\midspace{10cm} \caption{Figure 4}   

Figure 4  indicates the phase portraits 
of $f_0$ and some of its perturbations.  Observe that a perturbation can create 
a new type of orbits which go through ``the gate'' between the two fixed points.
Such orbits can give rise to a drastic change of the global dynamics (such as 
the inflation of the Julia set).   However it takes an extremely long time 
for these orbit to goes through the gate, so we need to consider a large number of 
iterates of the map in order to see the phenomenon.  

\smallskip 

  We analyze such a bifurcation using the theory of Ecalle cylinders, and 
the principle can be summarized as follows.  
For $f_0$, one can find ``fundamental regions'' $S_0^+$ and $S_0^-$, 
each of which has a boundary consisting of 
two curves joining the fixed point $0$  
such that one is mapped to the other.  See Figure 5. 

\midspace{6cm} \caption{Figure 5}

Gluing these two curves of $S_0^-$ (resp. $S_0^+$), one 
obtains a topological cylinder $\Cal C_0^-$ (resp. $\Cal C_0^+$), called the 
{\it outgoing} (resp. {\it incoming}) {\it Ecalle cylinder}, 
which turns out to be conformally 
isomorphic to the bi-infinite cylinder $\C^*$ (or $\CZ$).  
The orbits going from the ends (``horns'') of $S_0^-$ to $S_0^+$ induce 
a continuous and analytic mapping $\cE_{f_0}$ from a neighborhood of 
the ends of $\Cal C_0^-$ to $\Cal C_0^+$ ({\it the Ecalle transformation}).  
The identification 
$\CZ \to \Cal C_0^-$ can be lifted and extended to a map 
$\phi_0$, defined on a subset of $\C$ 
(which is considered to be the universal cover of $\CZ$) into the dynamical plane 
of $f_0$, 
similarly the identification  $\Cal C_0^+ \to \CZ$  can be lifted to a map 
$\Phi_0:\Cal B \to \C$, 
where $\Cal B$ is the parabolic basin of $0$.  Note that these functions can have 
critical points after the extension to the maximal domain of definition.  

\smallskip 

We consider the perturbation of the form $f(z)=e^{2 \pi i \alpha}z+O(z^2)$  with 
$\alpha \neq 0$ satisfying $|\arg \alpha|<\pi/4$.  
Then it can be shown that fundamental regions $S_f^-$, $S_f^+$ continue 
to exist, with boundary curves joining two fixed points.  
See Figure 5.  The quotient 
cylinders $\Cal C_f^-$, $\Cal C_f^+$ are still isomorphic to 
$\C^*$.  So we can define functions $\phi_f$, $\Phi_f$, $\cE_f$ which are similar 
to $\phi_0$, $\Phi_0$, $\cE_{f_0}$.  (The domains of definition may be smaller.)  

Now there is ``a gate'' open between the fixed points, 
and any orbit starting from $S_f^+$ passes through the gate and 
eventually falls into $S_f^-$.   
This induces a new map $\chi_f$ from $\Cal C_f^+$ to $\Cal C_f^-$, which 
is a conformal isomorphism.  
Thus we define {\it the return map} $\cR_f=\chi_f \circ \cE_f$, which 
corresponds near the ends of $\Cal C_f^-$ to the map of $S_f^-$ 
sending $z\in S_f^-$ to the first return point to $S_f^-$ along its forward 
orbit.   Therefore the return map corresponds to a heigh iterate of the map $f$.  
So one can study orbits of 
$f$ which return many times to the neighborhood of $0$ using  
the return map.  

  Moreover it will be shown that when $f$ tends to $f_0$ with the above restriction 
on $\arg \alpha$, 
the limit behaviour of the return map is determined by $\alpha$ and $\cE_{f_0}$, 

\medskip

  The Ecalle cylinders were first studied and applied to some problems by 
Douady-Hubbard and Lavaurs ([DH], [L]).  
The aim of this section is to state some notions and facts in this theory,  
in fact we state only the facts about $\phi_0$, $\Phi_0$, $\cE_{f_0}$, 
$\phi_f$ and $\cR_f$.   
The formulation presented 
here is somewhat different from [DH], [L].  
In this paper, we focus more on the return map, 
the quantitative aspect of the theory 
and the renormalization arising from the Ecalle transformation.  
The proof of these facts will be given in the Appendix, 
although most of the results can be found in [Mi] and [DH].

\bigskip
\subheading{Notations}

In the following, a function $f$ is always associated with its domain of 
definition $\Dom(f)$ (an open set of $\CC$ or $\CZ$ etc.); that is, two functions 
are considered as distinct, if they have distinct domains, even if one is an 
extension of the other.  
A neighborhood of an analytic map 
$f:\Dom(f) \to \CC$ is a set containing 
$$\{ g:\Dom(g) \to \CC| \ g \text{ is analytic}, \Dom(g) \supset K \  
\text{ and } \sup_{z \in K} d(g(z),f(z)) < \epsilon \}, $$ 
where $K$ is a compact set in $\Dom(f)$, $\epsilon > 0$ and $d( \cdot, \cdot)$ is 
the spherical metric.  (If the map is to some other space, then $d$ should 
be replaced by an appropriate metric.)  
The system of these neighborhoods defines ``the compact-open topology together 
with the domain of definition'',  which is unfortunately 
not Hausdorff, since an extension of $f$ is 
contained in any neighborhood of $f$.  

Let 
$$\Cal F=\{ f:\Dom(f) \to \CC |
\ f \text{ is analytic },  0 \in \Dom(f) \subset \CC \text{ and } 
f(0)=0 \}.  $$

We use the following notations:  

$\pi:\C \to \C^* \equiv \C -\{0\}$,  $\pi(z)=\exp(2 \pi i z)$; 

$\pi_1:\C \to \CZ$ (the natural projection); \ 
$\pi_2:\CZ \to \C^*$,  $\pi=\pi_2 \circ \pi_1$; 
\newline 
Note that $\pi_2$ sends ``the upper end'' of $\CZ$ ($\Im z \to +\infty$) to $0$, 
and ``the lower end'' ($\Im z \to -\infty$) to $\infty$.  

$T:\C \to \C$,  $T(z)=z+1$; 

$\t_0(z)= -\frac 1z$.

\bigskip

\subheading{4.1 Parabolic fixed point}  

Let  $f_0 \in \Cal F$ such that $f'_0(0)=1$ and $f''_0(0) \neq 0$.  
By a linear coordinate change, we may suppose that $f''_0(0)=1$.  
We assume this throughout in this section.  

\smallskip 

\proclaim\nofrills{Fact: } \  
There exist objects $\Omega_\pm, \phi_0, \Phi_0$ and $\cE_{f_0}$ 
satisfying (4.1.1)-(4.1.4).  
\endproclaim 
See Figure 6.  

\midspace{8cm} \caption{Figure 6}

\medskip 
\noindent 
(4.1.1) $\Omega_+$ and $\Omega_-$ are simply connected domains and 
$\overline \Omega_+, \overline \Omega_- \subset \Dom(f_0)$; 

the boundaries $\partial \Omega_+$, $\partial \Omega_-$ are Jordan curves 
containing $0$; 

$\Omega_+ \cup \Omega_- \cup \{0\}$ is a neighborhood of $0$; 

$f_0(\overline \Omega_+) \subset \Omega_+ \cup \{0\}$ 
\  and \  
$f_0(\Omega_- \cup \{0\}) \supset \overline \Omega_- $ ; 
 
$f$ is injective on $\Omega_+ \cup \Omega_- \cup \{0\}$; 
 
$\Omega_+ \cap \Omega_-$ consists of two components; 
 
$f_0^n \to 0$ as $n \to \infty$ uniformly on $\Omega_+$; 

a point $z$ belongs to the parabolic basin $\Cal B$ of $0$ for $f_0$, 
if and only if for some $n \geq 0$, $f_0^n(z)$ is defined and 
belongs to $\Omega_+$.

\smallskip 

Here, the {\it  parabolic basin} of a parabolic fixed point $\zeta$ 
of an analytic function $f$ is 
$$
\Cal B=\biggl\{ z \biggm| 
\matrix \ z \text{ has a neighborhood on which } 
f^n\ (n=1,2, \dots) \text{ are defined} \\
\text{ and } f^n \to \zeta \  
\text{ uniformly as }  n\to \infty 
\endmatrix \biggr\}.  
$$ 
Note that $\zeta$ itself is not in the parabolic basin.

\smallskip 
\noindent 
(4.1.2) $\phi_0: \Dom(\phi_0) \to \CC$ is an analytic function satisfying: 
$$
\phi_0(w+1)=f_0 \circ \phi_0(w) \ \ \text{ if both sides are defined}, 
$$ 
and in fact the left hand side is defined if and only if the right hand side is; 
\noindent 
$\Dom(\phi_0)$ contains 
$\Q_0=\{ w \in \C| \  \pi/3< \arg(w+\xi_0) < 5\pi/3 \ \} $ 
and $\{w| \ |\Im w|> \eta_0 \}$  for large $\xi_0, \eta_0 >0$;  

\noindent 
$\phi_0(\Q_0)=\Omega_-$; \  
$\phi_0$ is injective on $\Q_0$;  

\noindent 
let $\phi^*_0=\t_0^{-1} \circ \phi_0$, then 
$$\phi^*_0(w)=w+a \log w+b + o(1) \ \text{ as } \Q_0 \ni w \to \infty,$$
where $a$, $b$ are some constants.  
 
It follows from the above condition that 
if $f_0$ is a rational map or an entire function, 
then $\Dom(\phi_0)=\C$.  

\medskip 
\noindent 
(4.1.3) $\Phi_0:\Cal B \to \C$ is an analytic function satisfying  
$$
\Phi_0 \circ f_0(z) = \Phi_0(z)+1 \ \text{ for } z \in \Cal B ;
$$ 
and $\Phi_0$ is injective on $\Omega_+\ (\subset \Cal B)$.   
  
\medskip 
\noindent 
(4.1.4) Let $\tilde \Cal B = \phi_0^{-1}(\Cal B)$, 
then $T(\tilde \Cal B) = \tilde \Cal B$ and $\tilde \Cal B$ contains 
$\{w \bigm| |\Im w| > \eta_0\}$ for some $\eta_0>0$.  

Now define 
$\tE_{f_0}: \tilde \Cal B \to \C$ by 
$$\tE_{f_0}=\Phi_0 \circ \phi_0.  $$
It satisfies 
$$\tE_{f_0}(w+1) = \tE_{f_0}(w)+1 \ \text{ for } w \in \tilde \Cal B.$$  
Hence $\cE_{f_0}=\pi \circ \tE_{f_0} \circ \pi^{-1}: 
\pi(\tilde \Cal B) \to \C^*$ is well-defined.   
Moreover it extends to $0$ and $\infty$ analytically by $\cE_{f_0}(0)=0$ 
and $\cE_{f_0}(\infty)=\infty$, and $\cE'_{f_0}(0) \neq 0$, 
$\cE'_{f_0}(\infty) \neq 0$.  
So $\Dom(\cE_{f_0})=\pi(\tilde \Cal B) \cup \{0, \infty\}$.  

\smallskip 

The map $\cE_{f_0}$ is called the {\it Ecalle transformation}, 
or the {\it horn map} (since it is defined near the ends of the cylinder).   
  
\medskip 
\noindent 
(4.1.5) Normalization.  
Note that $\phi_0(w+c)$ and $\Phi_0(z)+c'$ also satisfy (4.1.2) and 
(4.1.3).  So we adjust $\Phi_0$ by adding a constant so that 
$$\cE'_{f_0}(0) = 1.$$

\bigskip 

\subheading{4.2 Perturbation}
Let 
$$\Cal F_1= \{ f \in \Cal F| \ f'(0)=\exp(2\pi i \alpha) 
\text{ with } \alpha \neq 0 \text{ and } |\arg \alpha|< \pi/4 \}.$$  
In the following, we consider only the perturbations $f \in \Cal F_1$. 
\definition{Notations}
If $f \in \Cal F$ and $f'(0) \neq 0$, we express the derivative $f'(0)$ as 
$$
f'(0)=\exp (2\pi i \alpha(f))
$$ 
where $\alpha(f) \in \C$ and $-\frac 12 <\Re \alpha(f) \leq \frac 12$.  

If $f\in \Cal F$ is close to $f_0$,  
then $f$ has two fixed points near $0$ counted with multiplicity, 
one of them is $0$ and the other is denoted by $\s(f)$.    Note that 
$$\s(f)=-2 \pi i \alpha(f)(1+o(1)) \text{ as } f \to f_0.$$    
\enddefinition 

\proclaim\nofrills{Fact: }  \ 
There exist a neighborhood $\N_0$ of $f_0$ and a constant 
$\xi_0>0$ (which can be the sama as that in the definition of $\Q_0$ (4.1.2)) 
such that for $f \in \N_0 \cap \Cal F_1$, there exist analytic maps 
$\phi_f:\Dom (\phi_f) \to \C$ and $\cR_f:\Dom(\cR_f) \to \CC$ 
satisfying (4.2.1)-(4.2.4).  
\endproclaim 
 
\smallskip 
\noindent   
(4.2.1) If $w, w+1 \in \Dom(\phi_f)$, then $\phi_f(w) \in \Dom(f)$ and 
$$\phi_f(w+1)=f \circ \phi_f(w) ; $$
$\Dom(\phi_f)$ contains $\Q_f$ (see Figure 6), where 
$$
\Q_f=\{ w \in \C \bigl|\ \pi/3< \arg(w+\xi_0) < 5\pi/3 \text{ and }
|\arg(w+ \frac 1{\alpha(f)}-\xi_0)|< 2 \pi/3 \ \};    
$$
\noindent 
$\phi_f(w) \to 0$ when $w \in \Q_f$ and  $\Im w \to +\infty$; \  
$\phi_f(w) \to \s(f)$ when $w \in \Q_f$ and  $\Im w \to -\infty$.      
 
\medskip  
\noindent   
(4.2.2) $0, \infty \in \Dom(\cR_f)$,  $\cR_f(0)=0$, $\cR_f(\infty)=\infty$; 
\ $\cR_f(\Dom(\cR_f)\cap \C^*) \subset \C^*$; 

$\cR'_f(0)=\exp(-2 \pi i \dfrac 1{\alpha(f)})$,  \ hence \ 
$\alpha(\cR_f) \equiv -\dfrac 1{\alpha(f)}\ (\text{mod } \Z)$.  

\medskip  
\noindent   
(4.2.3) If $w, w' \in \Dom(\phi_f)$,  
$\cR_f(\pi(w))=\pi(w')$ and $|\arg(w'+\frac 1{2\alpha(f)}-w)|<2\pi/3$, 
then $f^n(\phi_f(w))=\phi_f(w')$ for some 
$n \geq 1$.  

\smallskip 

Moreover if $U, U'$ are connected subsets of $\Dom(\phi_f)$ such that 
\ 
$\cR_f^m(\pi(U)) \subset \pi(U')$ for some $m \geq 1$,  
\ 
$\phi_f|_U$, $\pi|_{U'}$ are injective, and  
\ 
$|\arg(w'+\frac 1{2\alpha(f)}-w)|<2\pi/3$ for $w \in U$, $w' \in U'$,  
\noindent 
then there exists an $n > m$ such that 
$$
f^n = \phi_f \circ (\pi|_{U'})^{-1} \circ \cR_f^m \circ \pi \circ (\phi_f|_U)^{-1} 
\ \text{  on }  \phi_f(U).  
$$   

\medskip 
\noindent 
(4.2.4) With respect to the topology defined in Notations in \S 4.0, we have 
$$
\phi_f \to \phi_0  \ \text{  and  } \ 
e^{2 \pi i/\alpha(f)} \cR_f \to \cE_{f_0}
\ \text{  when } f \in \N_0 \cap \Cal F_1 \text{ and } f \to f_0.  
$$ 
 
Denote $\hE_{f_0}=\pi_2^{-1} \circ \cE_{f_0} \circ \pi_2$  and  
$\hR_f=\pi_2^{-1} \circ \cR_f \circ \pi_2$.  Then 
$$\hR_f + \frac 1{\alpha(f)} \to \hE_{f_0} 
\ \text{  when } f \in \N_0 \cap \Cal F_1 \text{ and } f \to f_0.  $$

\medskip 

\noindent
(4.2.5) To study the quantitative aspect, 
it is often more convenient to work with $\phi^*_f=\t_0^{-1} \circ \phi_f$.  

For $\eta \in \R$, define $D(\eta)=\{ w \in \C \bigl|\ |w- i\eta|\leq |\eta|/4 \}$ 
and 
$D'(\eta)=\{ w \in \C \bigl|\ |w- i\eta|\leq |\eta|/2 \}$.   

As a corollary of the above facts, we have: 
\medskip 
\noindent  
For large $\eta >0$, 
there exists a neighborhood 
$\N_1(\eta) (\subset \N_0)$ of $f_0$, depending on $\eta$, 
such that if $f \in \N_1(\eta) \cap \Cal F_1$,   
then $\phi_f$ is defined and 
injective on $D(\pm \eta)$, the image $\phi^*_f(D(\pm \eta))$ is 
contained in $D'(\pm \eta)$ and $\frac 12 <|(\phi^*_f)'|<2$ on $D(\pm \eta)$.

\bigskip

\subheading{4.3 Remarks} 

\noindent 
(4.3.1) For $f \in \Cal F$ with $3\pi/4<\arg\alpha(f)<5\pi/4$, instead of 
being in $\Cal F_1$, we can obtain a similar result with the following changes: 

the lower end of $\CZ$ corresponds to the fixed point $0$; 

in (4.2.1), $\phi_f(w) \to \s(f)$ ($w \in \Q_f$ and $\Im w \to +\infty$),  
$\phi_f(w) \to 0$ ($w \in \Q_f$ and  $\Im w \to -\infty$);  

in (4.2.1), (4.2.3) and (4.2.4), use $-\frac 1{\alpha(f)}$ 
instead of $\frac 1{\alpha(f)}$; 

in (4.2.2), $\cR'_f(\infty) = \exp(-2\pi i \frac 1{\alpha(f)})$.

\medskip 
\noindent 
(4.3.2) 
The above facts suggest the following 
factorization:  
$$\hR_f = 
\hat \Cal T_{-1/\alpha} \circ (\hR_f + \frac 1\alpha), 
$$
where $\hat \Cal T_{-1/\alpha}(w)=w-1/\alpha$ on $\CZ$, $\alpha=\alpha(f)$.  
Here $(\hR_f + \frac 1\alpha)$ is a non-linear map defined only 
in a subset of $\CZ$, but approaches a fixed map $\hE_{f_0}$,  
whereas, $\hat \Cal T_{-1/\alpha}$ is an isomorphism of $\CZ$, but 
depends sensitively on $f$, since $\alpha(f) \to 0$ as $f \to f_0$.  
Therefore if $\{f_n\}$ is a sequence in $\N_0 \cap \Cal F_1$ such that 
$f_n \to f_0$ and $1/\alpha(f_n)-k_n \to -c$ as $n \to \infty$, where 
$k_n$ are integers, 
then there exists a limit 
$$\lim_{n \to \infty} \hR_{f_n} =\hE_{f_0} +c,$$ 
since we are considering the maps in $\CZ$.  
\smallskip 

The Ecalle transformation $\cE_{f_0}$ was originally defined as a map 
between two different spaces--- (a part of) $\Cal C_0^-$ and $\Cal C_0^+$.  
So there is, a priori, no meaning as a dynamical system.  However 
the above observation implies that for any $c \in \C$, 
the map 
$$ w \longmapsto \hE_{f_0}(w)+c $$ 
on $\CZ$ can be realized as the limit of return maps of $f_n$.  

\medskip 
\noindent
(4.3.3) 
In [DH] and [L], they use $\phi_0 \circ \Cal T_c \circ \Phi_0$, 
where $\Cal T_c(z)=z+c$ $(c \in \C)$,  instead of 
$\tE_{f_0}+c=\Cal T_c \circ \Phi_0 \circ \phi_0$.  
Then for a sequence $\{f_n\}$ as above, there is a sequence of integers $k_n$ 
such that 
$\phi_0 \circ \Cal T_c \circ \Phi_0$ is the limit of $f_n^{k_n}$ on 
the parabolic basin $\Cal B$.

\medskip 
\noindent
(4.3.4) A more accurate statement of (4.2.3) is as follows: 
There exists a particular lift of $\cR_f$ to the universal cover, 
$\tR_f:\pi^{-1}(\Dom(\cR_f)) \to \C$, such that 
if $w, w' \in \Dom(\phi_f)$ and $\tR_f^m(w)+n=w'$, where $n, m\geq 0$ are integers, 
then $f^n(\phi_f(w))=\phi_f(w')$.

\bigskip 
\bigskip

\heading
\S 5. Global behaviour of Ecalle transformation
\endheading 

\medskip 

In this section, we study the property of the maximal extension of 
the Ecalle transformation for a certain class of maps $\Cal F_0$.  
As in the previous section, we suppose that $f_0 \in \Cal F$ is a function 
satisfying $f'_0(0)=1$ and $f''_0(0) = 1$.  So we have $\phi_0$ and $\cE_{f_0}$ 
as in the previous section.  Let us denote $g_0=\cE_{f_0}$.  

\subheading{Inverse orbits and $\phi_0$}  

\definition{Definition} 
For a mapping $f$, a sequence of points $\{z_j\}_{j=0}^\infty$ is called 
{\it an inverse orbit} (of $z_0$) for $f$, if $z_j \in \Dom(f)$ and 
$f(z_j)=z_{j-1}$ for $j\geq 1$.  
\enddefinition 

\proclaim{Lemma 5.1}
For $w \in \Dom(\phi_0)$, let $z_j=\phi_0(w-j)$ $(j=0,1,\dots )$.  
Then  $\{z_j\}_{j=0}^\infty$ is an inverse orbit for $f_0$ converging to $0$.  
This gives a one to one correspondence between $\Dom (\phi_0)$  
and the set of inverse orbits converging to $0$, except 
the orbit $z_j=0$ $(j=0,1,\dots)$.  

Moreover if $z_j\ (j \geq 1)$ are not critical points, then 
$\phi'_0(w) \neq 0$.   
\endproclaim 

\demo{Proof} If $w \in \Dom(\phi_0)$, then $w-j \in \Dom(\phi_0)$ $(j=0,1,\dots )$ 
and $\phi_0(w-j) \to 0$ by (4.1.2), so 
the first statement is obvious.  
Let $\{z_j\}$ be an inverse orbit converging to $0$ and suppose $z_j \neq 0$ 
for some $j$.  
For large $j$, say for $j \geq j_0$, 
$z_j$ belongs to $\Omega_+\cup \Omega_-$ (see (4.1.1)).  
But $f_0^n$ tends to $0$ uniformly on $\Omega_+$, so  
there exists $j_0$ such that $z_j \in \Omega_-$ for 
$j \geq j_0$.  Let $w=(\phi_0|_{\Q_0})^{-1}(z_j)+j$ ($j\geq j_0$).  
It is easy to see that $w$ does not depend 
on $j\geq j_0$ and corresponds to the inverse sequence $\{z_j\}$.  
If $w$ and $w'$ give the same inverse sequence,  
take $j\geq 0$ such that 
$w-j, w'-j \in \Q_0$, 
then we have $w=w'$ by the injectivity of $\phi_0$ on $\Q_0$.  

The last statement follows from the facts that 
$\phi_0(w)= f_0^n \circ \phi_0 (w-n)$ for $w \in \Dom(\phi_0)$ and 
that $\phi'_0 \neq 0$ in $\Q_0$.  
\qed \enddemo 

\medskip 

\subheading{Immediate parabolic basin} 

\definition{Definition} 
Let $\zeta$ be a parabolic fixed point of an analytic function $f$.  
A connected open set $B \subset \Dom(f)$ is called {\it an   
immediate parabolic basin} of $\zeta$ for $f$,  if 
$B$ is a connected component of the parabolic basin of $\zeta$ for $f$ 
(see (4.1.2)) such that 
$f(B) = B$  and  $f:B \to B$ is proper (hence a branched covering).  

For a rational map, a parabolic periodic point always has an 
immediate parabolic basin and it coincides with the previous definition 
(\S 1).  
However a parabolic fixed point in general may not have any immediate 
parabolic basin.  
\enddefinition 

\proclaim{Lemma 5.2} 
Suppose that $f \in \Cal F$ and $f'(0) =1$, $f^{(k)}(0)=0 \ (1<k\leq q)$, 
$f^{(q+1)}(0) \neq 0 \ (q \geq 1)$, i.e. $q+1$ is the order of $f(z)-z$ at $0$.  
Then $f$ has at most $q$ immediate parabolic basins of $0$, 
each of which contains at least one critical value of $f$, 
except for a parabolic M\"obius transformation.    
Moreover if an immediate parabolic basin contains only one critical point 
(or critical value), then it is simply connected.  
\endproclaim 

\demo{Proof} 
This is a well-known argument for rational maps.  Let us see it briefly.  
By local analysis as in \S 4.1 (case $q=1$) (see also \S 7 or [Mi]), 
it can be shown 
that there exist $q$  disjoint simply connected open sets 
$V^{(1)}, \dots , V^{(q)}$ $(\subset \Dom (f))$ (called ``attracting petals'') 
such that \ 
$0 \in \partial V^{(i)}$;\ $f(\overline V^{(i)}) \subset V^{(i)}\cup \{0\}$; \  
$f$ is injective on $V^{(i)}$; \ 
a point $z \in \Dom(f)$ belongs to 
the parabolic basin of $0$ if and only if $f^n(z) \in \cup_i V^{(i)}$ for some 
$n \geq 0$; and the orbit spaces $V^{(i)}/_\sim$ (where $z \sim f(z)$ 
if $z, f(z) \in V^{(i)}$), are 
isomorphic to $\C$.  

Let $B$ be an immediate basin of $0$, and $B'$ the parabolic basin.  
By the above, $B \cap \cup_i V^{(i)} \neq \emptyset$ and $\cup_i V^{(i)} \subset B'$.   
Since $B$ is a component of $B'$, $B$ contains one of $V^{(i)}$'s.  So there are 
at most $q$ immediate basins.  
 
Define $V^{(i)}_n$ $(n=0,1, \dots)$ inductively, as follows: 
$V^{(i)}_0=V^{(i)}$; if $f(V^{(i)}_n) \subset V^{(i)}_n$, $V^{(i)}_{n+1}$ is 
the component of $f^{-1}(V^{(i)}_n)$ containing $V^{(i)}_n$, which exists 
and $f(V^{(i)}_{n+1})=V^{(i)}_n \subset V^{(i)}_{n+1}$.  
It is easy to see that for each $i$, $\cup_{n \geq 0} V^{(i)}_n$ is a component 
of $B'$.  
Now suppose $B=\cup_{n \geq 0} V^{(i)}_n$ is an immediate parabolic basin, hence 
by definition, $f:B \to B$ is a branched covering.  If $B$ contains 
no critical value or no critical point, 
then $f: V^{(i)}_{n+1} \to V^{(i)}_n$, hence $f^n:V^{(i)}_n \to V^{(i)}$ 
are covering maps, therefore they induces a covering map $B \to V^{(i)}/_\sim$.  
Moreover $V^{(i)}_n$ hence $B$ are simply connected.  
Therefore $B$ must be isomorphic to 
$\C$ and $f$ is a M\"obius transformation.  

Similarly, if $B$ contains only one critical point or critical value, then 
one can show inductively that $V^{(i)}_n$ are simply connected, hence so is $B$.  
\qed \enddemo

Let us go back to our $f_0$ as at the beginning of this section.  
Suppose that $f_0$ has an immediate parabolic basin $B$ of $0$.  
The above argument applies to $V^{(1)}=\Omega_+$, 
hence $B$ is unique and contains $\Omega_+$ (see (4.1.1)).  
Let $B'$ be the (whole) parabolic basin of $0$, 
and define $\tilde B= \phi_0^{-1}(B)$ and $\tilde B'=\phi_0^{-1}(B')$.  
By (4.1.4), $\{w \bigl|\ |\Im w| > \eta_0 \} \subset \tilde B'. $
Since for any $w \in \tilde B'$ 
there is $n \geq 0$ such that $T^n(w) \in \tilde B$, we have 
$$\{w \bigl|\ |\Im w| > \eta_0 \} \subset \tilde B. $$ 

Denote by $\tilde B^u$ (resp. by $\tilde B^\ell$) the component of $\tilde B$ 
containing $\{ w |\ \Im w > \eta\}$ (resp. 
$\{ w |\ \Im w < -\eta\}$).  
(The superscript ``$u$'' stands for upper, and ``$\ell$'' for lower.)  
In general,  $\tilde B^u$ and $\tilde B^\ell$ may coincide.  
Obviously $T\tilde B^u=\tilde B^u$, $T\tilde B^\ell=\tilde B^\ell$.  
Then define $B^u = \pi(\tilde B^u)$, $B^\ell = \pi(\tilde B^\ell)$.  

\bigskip 

\subheading{Covering property of $g_0$ }  

\definition{Definition} 
Let $\Cal F_0$ be the set of functions $f \in \Cal F$ such that 
$f'(0)=1,\ f''(0)=1$  and  $f$ has an immediate parabolic basin 
which contains only one critical point of $f$.  
Then it is automatically simply connected by Lemma 5.2.  
\enddefinition 

\proclaim{Proposition 5.3} 
Let  $f_0 \in \Cal F_0$ and $B, \tilde B^u, \tilde B^\ell$ as before.  Then    
$g_0:B^u \cup B^\ell \to \C^*$ is a branched covering of 
infinite degree, ramified only over one point $v \in \C^*$.   

The sets $\tilde B^u$, $\tilde B^\ell$, $B^u \cup \{0\}$, 
$B^\ell \cup \{\infty \}$ 
are simply connected, and  $B^u \cap B^\ell = \emptyset$.  
\endproclaim 

\demo{Proof} 
Let us first show that $\Phi_0: B \to \C$ is a branched covering (see (4.3.2)).  
In fact, $\Phi_0|_{\Omega_+}$ is injective, hence 
$\Phi_0:f_0^{-n}(\Omega_+)\cap B \to T^{-n}\Phi_0(\Omega_+)$ is a branched covering 
$(n\geq 0)$.  So it follows that $\Phi_0$ on $B$ is a branched covering.  

Now let us show that $\phi_0:\tilde B^u \cup \tilde B^\ell \to B$ is 
a branched covering.  
Let $z$ be a point in $B$.  Take simply connected neighborhoods $U$, $U'$ of $z$ 
such that $\overline U \subset U'$ and $U'$ contains at most one 
forward orbit of the critical point.  Let 
$z'$ be a point in $\phi_0^{-1}(U) \cap (\tilde B^u \cup \tilde B^\ell)$.  
Let $U'_n$ be the component of $f_0^{-n}(U')$ containing $\phi_0(z'-n)$ 
$(n=0,1,\dots)$.  Then for some $m \geq 0$, 
$U'_n$ $(n \geq m)$ do not contain the critical point.  Hence
there exist inverse branches $f_{0, U'_m}^{(-k)}:U'_m \to U'_{m+k}$ of $f_0^k$.  
The family $\{f_{0, U'_m}^{(-k)}\}$ is normal, since it avoids 
at least three values ($0$ and the orbit of the critical point).  
Moreover it converges 
to $0$ uniformly on compact sets, since it does so near $\phi_0(z'-m)$.  
Hence there exists an $n \geq m$ such that $U_n=f_{0, U'_m}^{(-n+m)}(U_m) 
\subset \Omega_-=\phi_0(\Q_0)$.  Let $V=T^n \circ (\phi_0|_{\Q_0})^{-1}(U_n)$.  Then 
$z' \in V$ and $\phi_0|_V=(f_0^n|_{U_n}) \circ (\phi_0|_{\Q_0}) \circ T^{-n}$.  
So $\phi_0:V \to U$ is a branched covering with at most one critical point, 
since $U$ contains at most one critical orbit.  
This shows that each component of $\phi_0^{-1}(U)$ is either unramified or 
ramified over a common point in $U$, 
hence  $\phi_0:\tilde B^u \cup \tilde B^\ell \to B$ is a branched covering.  

Therefore $\tE_{f_0}=\Phi_0 \circ \phi_0:\tilde B^u \cup \tilde B^\ell \to \C$ 
and $\cE_{f_0}:B^u \cup B^\ell \to \C^*$ are branched coverings.  
Moreover it is easy to see from the above that $\cE_{f_0}$ is ramified 
only over $v=\pi \circ \Phi_0(c)$, where $c$ is the unique critical point.  

Now let $U$ be a component of $f_0^{-1}(f_0(\Omega_+))$ contained in $B$, different 
from $\Omega_+$.  Then we have $f_0^n(U) \cap U = \emptyset$ $(n \geq 1)$.  
Hence for any component $V$ of $\phi_0^{-1}(U)$, $\pi|_V$ is injective (by the 
functional equation for $\phi_0$).  On the other hand, 
$\pi \circ \Phi_0:U \to \C^*$ is infinite to one, since 
$\pi \circ \Phi_0|_U = \pi \circ \Phi_0|_{\Omega_+} \circ f_0|_U$ and 
$\pi \circ \Phi_0|_{\Omega_+}$ is infinite to one.  
Hence $g_0$ is of infinite degree.  

\smallskip 

Let us show the simple connectivity of $\tilde B^u$ (or $\tilde B^\ell$).  
If $\gamma$ is a closed curve in $\tilde B^u$, 
$\gamma'=T^{-n}\gamma \subset \Q_0$ for some $n \geq 0$.  
Let $W$ be a region bounded by $\phi_0(\gamma')$, not containing $0$.  
Since both the immediate basin $B$ and $\Omega_-=\phi_0(\Q_0)$ are simply connected 
and do not contain $0$, we have $W \subset B \cap \Omega_-$.  
It follows that  $\gamma'$ is trivial in $\tilde B^u$ and so is $\gamma$.   
Therefore $\tilde B^u$ is simply connected.  

Then $B^u \cup \{0\}$ (or $B^\ell \cup \{0\}$) is also simply connected, since 
the fundamental group of $B^u$ is generated 
by a curve around $0$, which is trivial in $B^u \cup \{0\}$.  

Finally, let us show that 
$\tilde B^u$ and $\tilde B^\ell$ are different components of $\tilde B$.   
Suppose $\tilde B^u=\tilde B^\ell$.  Then $\tilde B^u=\tilde B^\ell=\C$,  
since $\tilde B^u$ and $\tilde B^\ell$ are invariant under $T$, simply connected, 
and contain half planes.  Therefore $B$ contains $\Omega_+$ and 
$\Omega_-=\phi_0(\Q_0)$, and the union is a punctured neighborhood of $0$.  
But $B$ is simply connecteda punctured neighborhood of $0$, so $B=\CC-\{0\}$.  
Hence $f_0$ is a parabolic M\"obius transformation, since it is analytic 
on $\CC$ and has no periodic point in $B$.  Then $f_0$ has no critical point 
and this contradicts with the assumption.  
Thus we have $\tilde B^u \cap \tilde B^\ell = \emptyset$,   
hence $B^u \cap B^\ell = \emptyset$.  
\qed \enddemo 

\medskip 

\subheading{Iteration of $g_0$}  

By the normalization in (4.2.5), we have $g'_0(0)=1$.  So $0$ is again 
a parabolic fixed point of $g_0$.  

\proclaim{Lemma 5.4} 
Let $f_0 \in \Cal F_0$ and $g_0=\cE_{f_0}$.  Then $g''_0(0) \neq 0$ and 
$g_0$ has a simply connected immediate parabolic basin which contains only one 
critical point.  In other words, $g_0$ belongs to $\Cal F_0$ after 
a linear scaling of the coordinate.  
Moreover $g_0^n(v)$ $(n=0,1,\dots)$ are 
defined and $g_0^n(v) \to 0$ $(n \to \infty)$, where $v$ is the unique 
critical value of $g_0$.  
\endproclaim 

\demo{Proof} Obviously $g_0 \not \equiv 0$.  
So there exist attracting petals $V^{(i)}$ $(i=1, \dots, q)$ for $g_0$ as in 
the proof of Lemma 5.2, where $q+1=ord (g_0(z) -z)$.   
Since $g_0:B^u \cup B^\ell \to \C^*$ is a branched covering, 
we can also construct $V^{(i)}_n$ and $g_0:V^{(i)}_{n+1} \to V^{(i)}_n$ 
is a branched covering with at most one critical point.  
Then 
$V^{(i)}_n$ are simply connected and $\deg g_0|_{V^{(i)}_n}$ $(n=0,1, \dots)$ 
are eventually constant.  It follows that $B_i=\cup_i V^{(i)}_n$ are simply 
connected and $g_0: B_i \to B_i$ is a branched covering with at most one 
critical point.  Hence $B_i$ are immediate parabolic basins.  
By Lemma 5.2, we have $q=1$, i.e., $g''_0(0) \neq 0$.  
The rest follows easily.  
\qed \enddemo

\bigskip 
\bigskip

\heading
\S 6. The construction of a hyperbolic subset
\endheading
\medskip 

In this section, we prove Theorem 2 in the case the multiplier is $1$  (hence 
$q=1$, $p=0$).  We will see in \S 7 how to modify the proof in other cases.  

To construct a hyperbolic subset as in \S 2, we trace certain inverse images 
of the fixed point $0$ using 
$\phi_0$, $\cE_{f_0}$, $\psi_0$, $\cE_{g_0}$, etc., then analyze the perturbed maps  
along these orbits.

\subheading{Step 0.  $f_0$ and $\{z_j\}$ }  

Suppose $f_0$ is a rational map and  $f_0 \in \Cal F_0$.   
There exists an inverse orbit (see \S 5 for definition) 
$\{z_j\}_{j=0}^{\infty}$ for $f_0$ such that $z_0=0$, $z_j \neq 0 \ (j >0)$ and  
$z_j \to 0 \ (j \to \infty)$.  In fact, since $0$ is in the Julia set, 
there is an inverse image of $0$ near $0$ (see [Bl], [Mi]).  
This point must  have an inverse 
sequence converging to $0$, otherwise it would belong to the parabolic basin.   
(see the local analysis in \S 4.1.)

\subheading{Step 1. $g_0$ and $\{w_j\}$, $\{w'_j\}$ }

By the construction in the previous sections, we obtain $\phi_0$ and 
$g_0=\cE_{f_0}$.  
Let $v$ be the unique critical value of $g_0$.  

By Lemma 5.1, there exists a unique $\tilde w_0 \in \Dom(\phi_0)$ 
corresponding to the inverse orbit $\{z_j\}$  such that 
$\phi_0(\tilde w_0-j)=z_j$.  In fact, for 
Moreover note that $z_0=0$  does not belong to 
the parabolic basin 
of $0$, hence $\pi(\tilde w_0) \notin \Dom(g_0)$ by the definition of $g_0$.  

\proclaim{Lemma 6.1} 
For given $w_0 \in \C^*$, there exist two inverse orbits 
$\{w_j\}_{j=0}^{\infty}$, $\{w'_j\}_{j=0}^{\infty}$  of $w_0$ for the map 
$g_0$  
such that: 
\newline 
$w_0=w'_0$;\ \ $w_j,\ w'_j \to 0  \ (j \to \infty)$; \  
$\{w_j\}_{j=1}^{\infty} \cap \{w'_j\}_{j=1}^{\infty} = \emptyset$; \    
and $w_j$, $w'_j$ $(j \geq 2)$ are not critical point of $g_0$.  
Moreover if $w_0 \notin \Dom(g_0)$, $w_1$, $w'_1$ are not critical point.  
\endproclaim

\demo{Proof}  
First note that $g_0^{-1}(w_0)$ can contain at most one periodic point.  
It is easy to see that if $g_0^{-1}(w_0)$ contains two points of 
$\{g_0^n(v)| n \geq 0 \}$, then one of them must be periodic.   Hence 
$\sharp (g_0^{-1}(w_0) \cap \{g_0^n(v)| n \geq 0 \}) \leq 2$.  
So we can choose two distinct non-periodic points $w_1$ and $w'_1$ in 
$(g_0)^{-1}(w_0) - \{g_0^n(v)| n \geq 0 \}$, 
since $g_0$ is a branched covering of infinite degree.  

Since $0$ is a parabolic fixed point 
of $g_0$, there exists an inverse orbit $\{w''_j\}_{j=1}^\infty$ 
for $g_0$  such that 
$w''_j \to 0$ as $j \to 0$, $w''_j \neq 0$  and  
$\{w''_j\}_{j=1}^\infty \cap \{(g_0)^n(v)\}_{n=0}^\infty =\emptyset$.  
By Lemma 5.4, $g_0^n(v) \to 0$ ($n \to \infty$).  
So we can take a simply connected open set 
$D \subset \C^*-\{g_0^n(v)\}_{n=0}^\infty $  containing $w_1$, $w'_1$ and 
$w''_1$.  

Since $g_0$ is a branched covering onto $\C^*$ ramified only over $v$, for each 
$j \geq 1$, there exists an inverse branch 
$G_j:D \to \C^*$  of $g_0^j$ such that 
$g_0^j \circ G_j = id_D$  and  $G_j(w''_1)=w''_{j+1}$.  
The family $\{G_j\}$ is 
normal, since it omits at least three values $0$, $\infty$ and $v$.  
As $w''_{j+1}=G_j(w''_1) \to 0$ ($j \to \infty$),  
$G_j \to \infty$.   

Now define $w_{j+1}=G_j(w_1)$ and $w'_{j+1}=G_j(w'_1)$.  
It is easy to check that these sequences have the claimed properties.  
\qed  
\enddemo

Applying this lemma to $w_0=\pi(\tilde w_0)$, we obtain $\{w_j\}$, $\{w'_j\}$ 
as above.  

\medskip

\subheading{Step 2. $\hat h$ and $\{W_j\}$ }

As $g_0 \in \Cal F_0$,  one can apply the construction of \S\S 4 and 5 
to $g_0$.  
Denote $\psi_0=\phi_{0,g_0}$ (``$\phi_0$'' corresponding to $g_0$) and  
$h_0=\Cal E_{g_0}$, for simplicity.    
As in the previous step, by Lemma 5.1, 
there exist $\tilde \zeta_0, \tilde \zeta'_0 \in \Dom(\psi_0)$  
corresponding to $\{w_j\}$, $\{w'_j\}$, respectively.  
Hence $w_0=\psi_0(\tilde \zeta_0)$ and $w'_0=\psi_0(\tilde \zeta'_0)$.   Moreover 
$\psi'_0(\tilde \zeta_0) \neq 0$ and $\psi'_0(\tilde \zeta'_0) \neq 0$.  
The following is the key lemma in the proof of Theorem 2.  

\proclaim{Lemma 6.2} 
Let $b>0$.  
There exist: a neighborhood $\N_2$ of $h_0$, where we consider that 
$0, \infty \in \Dom(h_0)$; two disjoint discs 
$W, W' \subset \CZ$ containing $\hat \zeta_0 =\pi_1(\tilde \zeta_0)$,  
$\hat \zeta'_0=\pi_1(\tilde \zeta'_0)$ respectively; and positive 
constants $C_0, C_1, C'_1$,  with the following properties.  
\newline 
If $h_1 \in \N_2$, $\beta \in \C/Z$ with $|\Im \beta|\leq b$ and 
$$
\hat h(\zeta)= \pi_2^{-1} \circ h_1 \circ \pi_2 (\zeta) - \beta \ \ 
(\text{for } \zeta \in  \Dom (\hat h) \equiv \pi_2^{-1} (\Dom (h_1)) \subset \CZ), 
$$ 
then there exists a sequence of disjoint topological discs $W_j \subset \CZ$ 
satisfying:

$W_0=W$ or $W'$;  

$W_j \subset \Dom(\hat h)$, $\hat h(W_j)=W_{j-1}$ and 
$\hat h|_{W_j}$ is injective $(j \geq 1)$;  

for any $K>0$, $W_j \subset \{ \zeta \in \CZ \bigl|\ |\Im \zeta| > K \}$, 
for large $j$;  

\comment 
(if $W_j \subset  \{ \zeta \in \CZ \bigl|\ |\Im \zeta| >R_3 \}$, then )  
\endcomment 

$diam W_j <1/2$,\ \   $dist(W_j,W_{j+1}) < C_0$, and 
  
$C_1<|(\hat h^j)'|<C'_1$ on $W_j$, for $j \geq 0$.  
\endproclaim

\demo{Proof}  Let us first consider 
$h=e^{2\pi i \beta}h_0$ with $|\Im \beta|\leq b$, and 
$\hat h=\pi_2^{-1} \circ h \circ \pi_2$.   We denote by  
$\tilde B^u$ and $\tilde B^\ell$ as in \S 5, the upper and the 
lower domain of definition of 
$\tE_{g_0}$ (not for $\tilde g_0=\tE_{f_0}$!) and also define 
$B^u=\pi(\tilde B^u)$ and $B^\ell=\pi(\tilde B^\ell)$

Since $B^u$ and $B^\ell$ are disjoint, at least one of them, say $B^u$, 
does not contain the unique value of $h$.  Hence, by Lemma 5.3, 
the local inverse of $h$ 
near $0$ can be extended to $B^u \cup \{0\}$.  So we have an analytic 
function $H:B^u \cup \{0\} \to B^u \cup \{0\}$ such that 
$h \circ H = id$, $H(0)=0$ and $|H'(0)|<1$ by Schwarz' lemma.  
Then it is well known (see [Mi]) that there exists a linearizing coordinate 
$L(z)$ such that 
$L$ is conformal near $0$, $L(0)=0$, $L'(0)=1$ and 
$L \circ H(z) = H'(0)\cdot L(z)$ near $0$.  
Passing to $\CZ$ model, where we denote $\hat H= \pi_2^{-1} \circ H \circ \pi_2$, 
we obtain the following: 

\smallskip 
\noindent 
(a) There exist constants $y_0 \in \R$, $C''_1, C''_2 > 0$ and 
an analytic function $\hat L: Y=\{ \zeta \in \CZ \bigl|\ \Im \zeta > y_0\} 
\to \CZ$ such that in $Y$, $\hat H$ is defined, $0<\Im(\hat H(\zeta)-\zeta)<C''_1$, 
$\hat L \circ \hat H = \hat L + a$, where $a=\frac 1{2 \pi i} \log H'(0)$ and 
$\Im a >0$, and ${C''_2}^{-1} < |\hat L'|<C''_2$.  

\smallskip 

Note that in the case where $B^\ell$ does not contain the critical value, we obtain 
a similar result with $\Im (\cdot)$ replaced by $-\Im (\cdot)$.  

As for $\hat \zeta_0$ and $\hat \zeta'_0$,  at least one of them, 
say $\hat \zeta_0$, is not the critical value.  So pick 
$\hat \zeta_1 \in \hat h^{-1}(\hat \zeta_0) \cap B^u$ and 
let $\hat \zeta_j = H^{j-1}(\hat \zeta_1)$ $(j \geq 2)$.  
Then there exists $j_0 \geq 1$ such that 
$\zeta_{j_0} \in Y$.  

\smallskip 
\noindent 
(b) There exist $j_0 \geq 1$, a small disc neighborhood $W_0$ of $\hat \zeta_0$ 
(or $\hat \zeta'_0$) and its inverse image $W_1, \dots, W_{j_0}$ such that 
$\overline W_j$ $(j=0, \dots, j_0)$ are disjoint, are contained in $B^u$, 
and contain no critical points of $\hat h$; 
$\hat h$ maps $W_j$ onto $W_{j-1}$ bijectively $(j=1, \dots, j_0)$; 
$W_0 \cap Y = \emptyset$ and $W_{j_0} \subset Y$.  

\smallskip 

Note that the properties (a) and (b) are stable under a perturbation, 
i.e., (a) and (b) still hold for $h=e^{2\pi i \beta'}h_1$ with  
$h_1$ near $h_0$ and $\beta'$ near $\beta$, and the constants are uniform 
in the neighborhood.  Then, by the compactness of $\{ \beta \in \C/\Z \bigl| 
\ |\Im \beta| \leq b \ \}$, there exist a neighborhood $\N_2$ of $h_0$ 
disjoint discs $W$, $W'$ in $\CZ$ containing $\hat \zeta_0$, $\hat \zeta'_0$, 
respectively, such that for $h_1 \in \N_2$ and $\beta$ with $|\Im \beta| \leq b$, 
(a) and (b) hold with uniform constants and $W_0=W$ or $W'$.  
Moreover we may assume that there are also uniform estimates for 
$(\hat h^j|_{W_j})'$, $diam W_j$, $dist(W_j,W_{j+1}$ $(0 \leq j \leq j_0)$.  
Now define $W_j=H^{j-j_0}(W_{j_0})$ $(j> j_0)$.  
Then we have uniform estimates on $(\hat h|_{W_j})'$ etc., since 
$|(\hat h^{j-j_0})'(\zeta)|=|\hat L'(\zeta)|\cdot 
|\hat L'(\hat h^{j-j_0}(\zeta))|^{-1}$.  
The rest of statements can be checked easily.  
\qed \enddemo

\medskip

\subheading{Step 3. Many $U_i$'s }  

The results in \S4 applies both to $f_0$ and to $g_0$.  So let us denote 
the objects $\phi_f$, $\cR_f$, $\N_1(\eta)$, etc. in \S 4.2, by 
$\phi_f$, $\cR_f$, $\N_1(f_0,\eta)$ for $f_0$, $f$, and 
$\psi_g$, $\cR_g$, $\N_1(g_0,\eta)$ for $g_0$, $g$.

Let $W_j$ be the discs in Lemma 6.2.   
Then there exists a constant $\gamma>0$ (depending only on $C_0$) such that 
for large $\eta>0$, if  
$f \in \N_1(f_0,\eta) \cap \Cal F_1$ and $g \in N_1(g_0, e^{2\pi \eta}) \cap 
\Cal F_1$, then there exist 
disjoint topological discs $U_1, \dots , U_N$, where 
$N \geq \gamma  \eta (e^{2 \pi \eta})^2 $,  
satisfying: 
 
$\overline U_i \subset \phi_f (D(\eta)) \subset D'(\eta)$; 
  
$V_i = \pi \circ (\phi_f|_{D(\eta)})^{-1} (U_i) \subset \psi_g (D(\eta')) 
\subset D'(\eta')$, 
where $\eta'=\pm e^{2 \pi \eta}$;  
 
$W_{j(i)} = \pi_1 \circ (\psi_g|_{D(\eta')})^{-1} (V_i)$ 
for some $j(i) \in {\Bbb N}$.   

\demo{Proof} 
Let us use $\phi^*_f = \t_0^{-1} \circ \phi_f$, $\psi^*_f = \t_0^{-1} \circ \psi_f$ 
and $\pi^* =\t_0^{-1} \circ \pi$.  By Lemma 6.2, the $W_i$'s tend to 
the upper or lower end of $\CZ$.  Suppose they tend to the upper end, 
and let $\eta'=e^{2 \pi i \eta}$.  (In the other case we take 
$\eta'=-e^{2 \pi i \eta}$.)  
 
Then the disc $D(\eta')$ contains entirely at least $\gamma (\eta')^2$ 
components of $\pi_1^{-1}(W_i)$'s for some constant $\gamma>0$, 
since the area of $D(\eta')$ is $(\pi/16)(\eta')^2$.  
By (4.2.5), $\psi_g$ maps these components into $D'(\eta')$ injectively.  
The inverse image of $D'(\eta')$ by $\pi^*$ consists of components, 
each of which is contained in a ``box'' 
$$\{ z \in \C \bigm|\ n<\Re z<n+1, \ \eta+\frac 1{2\pi} \log \frac 12 
<\Im z < \eta+\frac 1{2\pi} \log \frac 32 \}$$
for some $n \in \Z$.  So $D(\eta)$ contains at least $\eta$ 
of these components, for $\eta$ large.  Finally $\phi^*_f$ maps 
$D(\eta)$ into $D'(\eta)$ injectively.  
As for the components of $\pi_1^{-1}(W_i)$, they have at least 
$\gamma \eta (\eta')^2$ entire preimages by 
$\pi_1 \circ (\psi^*_g|_{D(\eta')})^{-1} \circ 
\pi^* \circ (\phi_f|_{D(\eta)})^{-1}$.  And this proves the above assertion,  
using $(\psi^*_g|_{D(\eta')})^{-1} \circ \pi^*=(\psi_g|_{D(\eta')})^{-1} \circ 
\pi$.  
\qed \enddemo

Let us denote $U^*_i=\t_0^{-1}(U_i)$, $\tilde V_i= (\phi_f|_{D(\eta)})^{-1}(U_i)$, 
$V^*_i=\t_0^{-1}(V_i)$.  Then we have maps 
$$
U_i      @> \t_0^{-1} >> U^*_i @> (\phi^*_f|_{D(\eta)})^{-1} >> 
\tilde V_i @> \pi^* >> V^*_i @> \pi_1 \circ  (\psi^*_g|_{D(\eta')})^{-1} >> 
 W_{j(i)} 
$$
and the estimates on the derivatives 
$$\frac 12 < \bigl|[(\phi^*_f|_{D(\eta)})^{-1}]'\bigr|  <2, \  
\frac 12 < \bigl|[\pi_1 \circ (\psi^*_g|_{D(\eta')})^{-1}]'\bigr|  <2, \ 
\frac 14 \eta^2 <|(\t_0^{-1}|_{U_i})'|<\frac 94 \eta^2 
$$
$$
\text{ and } \pi e^{2 \pi \eta}< |(\pi^*|_{\tilde V_i})'| < 3\pi e^{2 \pi \eta}. $$ 
For the last estimate, use the fact that $\tilde V_i$ is contained in a box as above.

\subheading{Step 4. from $W_0$ to $U$ } 

There exists an integer $k \geq 0$ such that $z_j$ $(j >k)$ are not 
critical point of $f_0$.  Let $\pi^{(-1)}$ be the local inverse of $\pi$ near 
$w_0$ such that $\pi^{(-1)}(w_0)=\tilde w_0 -k$.  
Similarly let $\pi_1^{(-1)}$ be the local inverse of $\pi_1$ near 
$\zeta_0$ and $\zeta'_0$ such that $\pi_1^{(-1)}(\zeta_0)=\tilde \zeta_0$ 
and $\pi_1^{(-1)}(\zeta'_0)=\tilde \zeta'_0$.  

Note that in Lemma 6.2, 
we may take $W$, $W'$ smaller without changing 
the statement.  
By Lemma 5.1, we have 
$\phi'_0(\tilde w_0 - k) \neq 0$, $\psi'_0(\tilde \zeta_0) \neq 0$ and 
$\psi'_0(\tilde \zeta'_0) \neq 0$.  

\medskip 

Then it follows that  
we can change $W$, $W'$ smaller, so that there exist 
a small neighborhood $U$ of $z_k$, 
neighborhoods $\N_3(f_0)$, $\N_3(g_0)$ of $f_0$, $g_0$ and constants $C_2, C'_2$  
such that 
\newline 
for $f \in \N_3(f_0) \cap \Cal F_1$ and $g \in \N_3(f_0) \cap \Cal F_1$, \ 
$\phi_f \circ \pi^{(-1)} \circ \psi_g \circ \pi_1^{(-1)}$  
is defined and injective on $W$ and on $W'$, and their images cover 
$U$;  and moreover the derivative has a bound   
$C_2<|(\phi_f \circ \pi^{(-1)} \circ \psi_g \pi_1^{(-1)})'|<C'_2$  on $W \cup W'$.  

\medskip 

\subheading{Step 5.  last $k$ iterate }

Let $U$ be as above and $\nu=\deg_{z_k} f_0^k$.  
Then it can be easily seen that there exists  
a neighborhood $\N_4(f_0, \eta)$ of $f_0$, depending on large $\eta>0$, 
such that 
for $f \in \N_4(f_0, \eta)$, 
there exists an open set $\Cal U' \subset U$ such that 
$f^k:\Cal U' \to \t_0(D'(\eta))$ is bijective and 
$$C_3 \left(\frac1\eta \right)^{\frac {\nu-1}\nu} 
<|(f^k)'|<C'_3 \left(\frac1\eta \right)^{\frac {\nu-1}\nu} \ \
\text{ on } \Cal U', $$ 
where $C_3, C'_3>0$ are constants independent of $\eta$.

\medskip 

\subheading{Step 6.  a hyperbolic set $X_f$ } 

Let $f_0$ be a rational map belonging to the class $\Cal F_0$,  
$b > 0$ and $\eta > 0$ large.  
First note that if an analytic function is close to $f_0$, then 
it can be conjugated to a function in $\Cal F$ close to $f_0$ 
by a translation near $id$.  
So for Theorem 2, we only need to consider the functions in $\Cal F$.  

Suppose that $f \in \Cal F$ is close to $f_0$  and 
$$\alpha(f) 
= \cfrac 1 \\
\ a_1 - \cfrac 1 \\
\ a_2 + \beta \ \endcfrac $$ 
with large positive integers $a_1, a_2$ and $\beta \in \C$ 
satisfying $0 \leq \Re \beta < 1$, $|\Im \beta|\leq b$.   
Other cases with different signes in the expression of $\alpha(f)$ 
can be treated similarly, by using the complex conjugates or 
by reformulating the procedure for the lower end of $\CZ$ instead of the upper end.  
See Remark (4.3.1).  

If  $f$ is close to $f_0$ and $a_1$ is large,  then 
$|\arg \alpha(f)| < \pi/4$, $\cR_f$ is defined and 
$e^{2 \pi i/\alpha(f)}\cR_f$ is close to $g_0=\cE_{f_0}$ \ 
(or $\hR_f+\frac 1{\alpha(f)}$ is close to $\hat g_0=\hE_{f_0}$).  
Let us denote $g= \cR_f$.   
If, moreover, $a_2$ is sufficiently large, 
then $g$ itself is close to $g_0$, $|\arg \alpha(g)|= 
|\arg \frac 1{a_2+\beta}|<\pi/4$,   
$h=\cR_g$ exists and $h_1=e^{2\pi i/\alpha(g)} \cR_g= 
e^{2 \pi i \beta}h$ is close to $h_0=\cE_{g_0}$.  
(Note here that, by (4.2.2), $\alpha(g) \equiv -1/\alpha(f) 
\equiv 1/(a_2+\beta)\ (\text{mod }\Z)$.)  
Therefore, for $f$ close to $f_0$, $a_1, a_2$ large, we have 
$f \in \N_1(f_0, \eta) \cap N_3(f_0) \cap N_4(f_0, \eta)$, 
$g \in \N_1(g_0, |\eta'|) \cap N_3(g_0)$ and $h_1 \in \N_2$.   
In particular, we can apply Lemma 6.2 to $\hat h=\pi_2^{-1} \circ h \circ \pi_2$ 
to obtain $W_j$.  Then let $U_i$, $U$ be as in Steps 3 and 4.  
\bigskip 
\midspace{15cm} \caption{Figure 7}

Now let us consider the following sequence of maps: 

$$ \aligned 
U_i      @> \ \t_0^{-1}\  >> U^*_i @> (\phi^*_f|_{D(\eta)})^{-1} >> 
\tilde V_i @> \pi^* >> V^*_i @> \pi_1 \circ  (\psi^*_g|_{D(\eta')})^{-1} >> 
& W_{j(i)}  \\ 
&  \CD @VV \hat h^{j(i)} V \endCD \\ 
f^k(U) @< f^k << U  \subset  U' 
@< \ \ \ \ \phi_f \circ \pi^{(-1)} \ \ \  << V 
@< \ \ \ \ \ \psi_g \circ \pi_1^{(-1)} \ \ \ \ << & W_0  
\endaligned 
\tag{*} 
$$
\noindent 
where $V=\psi_g \circ \pi_1^{(-1)} (W_0) $ and 
$U'= \phi_f \circ \pi^{(-1)}(V)$.  
Note that all maps except the last $f^k$ are injective.  
Write $\Cal U=\t_0(D'(\eta))$ and let $\Cal U' \subset U$ be the set 
obtained in Step 5.  Then tracing the inverse image of $\Cal U$ via $\Cal U'$ 
by the maps (*), we obtain $\Cal U_i \subset U_i$.  

Moreover the composition of (*) on $\Cal U_i$ is equal to $f^{n_i}$ 
for some integer  $n_i \geq 1$ \ $(i=1, \dots, N)$.  In fact, 
by (4.2.3), $g^{m_i}=\psi_g \circ 
\pi_1^{(-1)} \circ \hat h^{j(i)} \circ 
\pi_1 \circ (\psi_g|_{D(\eta')})^{-1}$ on $V_i$, since 
$|\arg(z'+\frac 1{2\alpha(g)}-z)|<2\pi/3$ for $z \in D(\eta')$, 
$z' \in \pi_1^{(-1)}(W \cup W')$), 
if $|\eta'|$ is large.  Similarly, 
we have $f^{n'_i}=\phi_f \circ \pi^{(-1)} \circ g^{m_i} \circ 
\pi \circ (\phi_f|_{D(\eta)})^{-1} $ on $U_i$, for some $n'_i \geq 1$.  
And recall that $\phi^*_f=\t_0^{-1} \circ \phi_f$, $\pi^*=\t_0^{-1} \circ \pi$, 
etc.  

\medskip 

Thus we obtained $(f, \Cal U, \Cal U_i)$ as described in \S 2, i.e., 
$\overline {\Cal U}_i \subset \Cal U$ such that $f^{n_i}: \Cal U_i \to \Cal U$ 
is bijective ($i=1, \dots, N$).    Moreover, combining all the above estimates, 
we have 
$$N \geq \gamma\eta (e^{2 \pi \eta})^2 $$ 
and 
$$|(f^{n_i}|_{\Cal U_i})'| \leq C\eta^{1+1/\nu} e^{2 \pi \eta}, $$  
where $\gamma$ and $C$ are positive constants independent of $\eta$.  

Hence by Lemmas 2.1 and 2.2, we have a hyperbolic subset $X_f$ for $f$ and 
an estimate for the Hausdorff dimension 
$$
\delta = \HD X_f \geq 
\frac {\log \gamma + \log \eta +4 \pi \eta}{\log C + 
(1+\frac 1\nu)) \log \eta +2 \pi \eta} .  
$$  
The right hand side tends to $\infty$ as $\eta \to \infty$.   
Thus we have proved Theorem 2 in the case $f'_0(0)=1$.  
\qed 

\remark{Remarks}
(i) Note that the maps in (*) other than $f^k|_U$, $\pi^*$, $\t_0^{-1}|_{U_i}$ 
have bounded derivatives, and that the effect of $\pi^*$ 
dominates others.  
\newline 
(ii) The procedure $\Cal F_0 \ni f_0 \longmapsto g_0=\cE_{f_0} \in \Cal F_0$ 
can be considered as a renormalization.  It it related, as we have seen, 
to the return map of $f$ near $f_0$.  One can interpret that 
$\pi \circ \phi_0^{-1}$ is a correspondence between the phase spaces of 
$f_0$ and ``its renormalization'' $g_0$, and that it has an exponential effect, 
because $f$ moves points extremely slowly near $0$ and require 
a large number of iterate for the return map.  In this sense, the renormalization 
procedure is essential in the above estimate.  
\newline 
(iii) It will be instructive to make a ``caricature'' (proposed by Curt McMullen) 
to understand the situation.  For simplicity, assume $k=0$ $(\nu=1)$, and 
pretend that the maps in (*) other than $\t_0^{-1}$ and $\pi^*$ were affine.  
In particular, $\hat h$ is supposed to be a translation on $\CZ$.  So it 
produces a one dimensional array of discs $W_j$, then $\pi_1^{-1}$ unwraps them 
to a two dimensional array.  These discs are squeezed by $(\pi^*)^{-1}$, 
and finally inverted by $\t_0$.  One can show that an invariant set 
produced by this system has dimension two.  
\newline 
(iv)
Note that the map $f_0$ need not to be a rational map.  In fact, it is enough 
to assume that $f_0 \in \Cal F_0$ and $0$ has an inverse orbit $\{z_j\}$ 
as in Step 0.  
\endremark

\bigskip 
\bigskip

\heading
\S 7. Parabolic fixed points with multiplier $\neq 1$  
\endheading
\bigskip

Let us consider an analytic function $f_0(z)$ near $0$ such that 
$$f_0(0)=0, 
f'_0(0)=\exp(2 \pi ip/q), $$ 
where $p, q \in \Z$, $q > 1$ and $(p,q)=1$.   
Then it is known that 
if $f_0^q(z) \not \equiv z$, then it has an expansion of the form 
$$
f_0^q(z)=z+a_{\nu q+1}z^{\nu q+1} +O(z^{\nu q +2}), 
$$
where $\nu$ is a positive integer and $a_{\nu q +1} \neq 0$.  
In the following, we assume that $\nu=1$, in other words, 
$(f_0^q)^{(q+1)} \neq 0$.   Then, as before, we may assume that $a_{q+1}=1$.  
The dynamics of $f_0$ and its perturbation is 
shown in Figure 8.  

\midspace{6cm} \caption{Figure 8}

To analyze the bifurcation, we need to consider $q$ incoming and $q$ outgoing 
Ecalle cylinders $\Cal C_0^{k,+}, \Cal C_0^{k,-}$ $(k \in \Z/q\Z)$, see Figure 9.  
Now the Ecalle transformations map 
the upper end of $\Cal C_0^{k,-}$ to that of $\Cal C_0^{k,+}$, 
the lower end of $\Cal C_0^{k,-}$ to that of $\Cal C_0^{k-1,+}$.  

\midspace{6cm} \caption{Figure 9} 

Consequently, the statements of \S\S 4-6 should be changed as follows.  
We only remark the part which is to be changed.  

\medskip 
\subheading{Changes in \S 4} 

\medskip 
\subheading{4.1} 
(4.1.1) There $2q$ regions $\Omega_+^{(k)}, \Omega_-^{(k)}$ $(k\in \Z/q\Z)$ 
instead of two regions  $\Omega_+$, $\Omega_-$;  

$\cup_k \Omega_+^{(k)} \cup \Omega_-^{(k)} \cup \{0\}$ is a neighborhood of $0$, 
on which $f_0$ is injective; 

$f_0(\overline \Omega_+^{(k)}) \subset \Omega_+^{(k+p)} \cup \{0\}$ 
\  and \  
$f_0(\Omega_-^{(k)} \cup \{0\}) \supset \overline \Omega_-^{(k+p)} $ ; 
 
$\Omega_+^{(j)} \cap \Omega_-^{(k)}$ is non-empty and connected, if 
$j=k$ or $j=k-1$,  it is empty otherwise; 

$f_0^n \to 0$ as $n \to \infty$ uniformly on $\Omega_+^{(k)}$; 

The parabolic basins $\Cal B^{(k)}$ are defined to be 
$\cup_{n \geq 0} f^{-nq}(\Omega_+^{(k)})$, then 
a point belongs to $\cup _k \Cal B^{(k)}$ if and only if it has a neighborhood 
on which  
$f^n\ (n=1,2, \dots)$ are defined and $ f^{nq} \to 0$ uniformly as $n\to \infty$.  

Let us fix a $k \in \Z/q\Z$.  

\smallskip 
\noindent 
(4.1.2) $\phi_0$, $f_0$, $\Omega_-$ should be replaced by $\phi_0^{(k)}$, 
$f_0^q$, $\Omega_-^{(k)}$.  As for $\phi^*_0$, define 
$$\phi^{(k)*}_0(w)=-\frac 1{q(\phi_0^{(k)})^q}, $$
then 
$$\phi^{(k)*}_0(w)=w+O\left( w^{\frac {q-1}q}\right).  $$

\smallskip 
\noindent 
(4.1.3) $\Phi_0$, $\Cal B$, $f_0$, $\Omega_+$ are to be replaced by 
$\Phi_0^{(k)}$, $\Cal B^{(k)}$, $f_0^q$, $\Omega_+^{(k)}$.  

\noindent 
(4.1.4)  Let $\tilde \Cal B^{(k,u)} = (\phi_0^{(k)})^{-1}(\Cal B^{(k)})$ and 
$\tilde \Cal B^{(k,\ell)} = (\phi_0^{(k)})^{-1}(\Cal B^{(k-1)})$.   
then $\tilde \Cal B^{(k,u)}, \tilde \Cal B^{(k,\ell)}$ are invariant under 
$T$ and 
$\tilde \Cal B^{(k,u)} \supset \{w \bigm| \Im w > \eta_0\}$, 
$\tilde \Cal B^{(k,\ell)} \supset \{w \bigm| \Im w < -\eta_0\}$ for some $\eta_0>0$.  

Define $\tE_{f_0}^{(k,u)}: \tilde \Cal B^{(k,u)} \to \C$, 
$\tE_{f_0}^{(k,\ell)}: \tilde \Cal B^{(k,\ell)} \to \C$ by 
$$\tE_{f_0}^{(k,u)}=\Phi_0^{(k)} \circ \phi_0^{(k)} \ \text{ and } \ 
\tE_{f_0}^{(k,\ell)}=\Phi_0^{(k-1)} \circ \phi_0^{(k)}.  $$
They satisfy 
$$\tE_{f_0}^{(k,u)}(w+1) = \tE_{f_0}^{(k,u)}(w)+1 \ \text{ for } 
w \in \tilde \Cal B^{(k,u)}, \text{ etc}; $$  
$\cE_{f_0}^{(k,u)}=\pi \circ \tE_{f_0}^{(k,u)} \circ \pi^{-1}: 
\pi(\tilde \Cal B) \cup \{0\} \to \C$ is well-defined and analytic, 
and  ${\cE_{f_0}^{(k,u)}}'(0) \neq 0$.   Similarly, $\cE_{f_0}^{(k,\ell)}
=\pi \circ \tE_{f_0}^{(k,\ell)} \circ \pi^{-1}: 
\pi(\tilde \Cal B) \cup \{\infty\} \to \CC-\{0\}$ is analytic, 
and  ${\cE_{f_0}^{(k,\ell)}}'(\infty) \neq 0$.

\medskip 

\subheading{4.2 Perturbation}
$$\Cal F_1= \{ f \in \Cal F| \ f'(0)=\exp(2\pi i \frac {p+\alpha}q) 
\text{ with } \alpha \neq 0 \text{ and } |\arg \alpha|< \pi/4 \}.$$  

Set 
$$
f'(0)=\exp (2\pi i \frac {p+\alpha(f)}q).  
$$ 
The periodic points of $f$ of period $q$ near $0$ are labelled so that 
$$\s^{(k)}(f)=(-2 \pi i \alpha(f))^{1/q}e^{2\pi i k/q}(1+o(1)) 
\ \text{ as } f \to f_0,$$    
where  $|\arg (-2 \pi i \alpha(f))^{1/q}|< \pi/q$, then 
$f(\s^{(k)}(f))=\s^{(k+p)}(f)$.    

The functions $\phi_f$, $\cR_f$, $\cE_{f_0}$ are to be replaced by 
$\phi_f^{(k)}$, $\cR_f^{(k,u)}$, $\cE_{f_0}^{(k,u)}$.  
The functionnal equation for $\phi_f$ becomes 
$$\phi_f^{(k)}(w+1)=f^q \circ \phi_f^{(k)}(w). $$
In (4.2.3) and (4.3.4), $f^n$ is to be replaced by $f^{nq+r}$, where 
$r$ is an integer such that $0<r<n$ and $rp\equiv -1$ $(mod\ q)$.  
Finally 
$$
\phi^{(k)*}_f(w)=-\frac 1{q(\phi_f(w))^q}
$$ 
satisfies (4.2.5).  

\medskip 

\subheading{4.3 Remarks} 
(4.3.1) If $f\in \Cal F$ satisfies $3\pi/4<\arg\alpha(f)<5\pi/4$, then 
we use $\cR_f^{(k,\ell)}$, $\cE_{f_0}^{(k,\ell)}$ instead of 
$\cR_f^{(k,u)}$, $\cE_{f_0}^{(k,u)}$.  

\medskip 

\subheading{Changes in \S 5} 

Let $f_0$ be as in \S 4.  Then $f_0$ has $q$ parabolic basins $\Cal B^{(k)}$ 
$(k \in \Z/q\Z)$.  Let $\Cal F_0$ be the set of such functions $f_0 \in \Cal F$ 
having an immediate parabolic basin $B^{(k)}$ in each $\Cal B^{(k)}$, 
containing only one critical point of $f_0^q$.  Let $\tilde B^{(k,u)}$ 
be the component of $\tilde \Cal B^{(k,u)}$ containing $\{w\bigm| \Im w > \eta_0\}$, 
and $B^{(k,u)}=\pi(\tilde \Cal B^{(k,u)})$.  

Then $g_0=\cE_{f_0}^{(k,u)}:B^{(k,u)} \to \C^*$ is a branched covering 
of infinite degree, ramified only over one point.  (Proposition 5.3) 

Other statements are similar.   

\medskip 

\subheading{Changes in \S 6} 

Note that $g_0$ is in $\Cal F_0$ in the sense of \S 4, that is $g'_0(0)=1$.  
So we only need to change $\phi_f$ as above and $\t_0^{-1}|_{U_i}$ to 
$$z \mapsto -\frac 1{qz^q}. $$  
Hence in Step 3, the estimate on the derivative of $\t_0^{-1}$ should be replaced by 
$$C\eta^{q+1} < (\t_0^{-1}|_{U_i})' <C'\eta^{q+1}, $$
and in Step 6, the estimate on $(f^{n_i}|_{\Cal U_i})'$ becomes 
$$|(f^{n_i}|_{\Cal U_i})'| \leq C\eta^{q+1/\nu}e^{2\pi \eta}.$$  

This is enough to prove that $\HD X_f \to 2$ as $\eta \to \infty$.

\bigskip 
\bigskip

\heading 
Appendix.  Proof of the properties of Ecalle cylinders 
\endheading 
\medskip 

In this Appendix, we give the proof or comments for the facts which were stated in 
\S4 and \S7.  Note that the facts in \S 4.1 can be found in [Mi] and 
most of the facts in \S 4.2 can be found in [DH].  We prove the facts in \S 4.2 
by introducing a new coordinate, and it makes clearer the relationship between 
the return map and the renormalization of functions with irrationally 
indifferent fixed points (cf. Yoccoz [Y]).  

\subheading{A.1 Coordinate changes} 

  Let $f_0(z)=z+z^2+\dots$ as in \S 4.1.  Let us introduce a new coordinate $w$ 
by $z=-1/w=\t_0(w)$, then $f_0$ corresponds to the map $F_0$ of the form 
$$F_0(w)=w+1+O(1/w). \tag{A.1.1}$$  
For functions near $f_0$, we introduce a new coordinate by the following.   

\proclaim{Lemma A.1.2}
There exist a neighborhood $\N'$ of $f_0$ in $\Cal F$ and a neighbourhood 
$\Cal V$ of $0$ in $\C$  
such that if $f \in \N'$ then $\overline {\Cal V} \subset \Dom(f)$ and 
$f(z)$ can be expressed as 
$$ 
f(z)=z+z(z-\sigma)u(z), \tag{A.1.3} 
$$
where $\sigma=\s(f)$ is a point in $\Cal V$ and $u(z)=u_f(z)$ is a non-zero 
holomorphic function defined in a neighborhood of $\overline {\Cal V}$.  
Hence $0$ and $\sigma(f)$ are the only fixed points of $f$ in $\Cal V$.  
Moreover $\s(f_0)=0$, $u_{f_0}(z)=(f_0(z)-z)/z^2$, 
$e^{2 \pi i\alpha(f)}=f'(0)=1-\s(f)u_f(0)$, hence 
$$\s(f)= - 2\pi i \alpha(f) (1+o(1)) \ \ \text{as } 
f \to f_0 .  
\tag{A.1.4}$$ 
The correspondence  
$f \mapsto \s(f)$, $f \mapsto u_f(z)$ (with $\Dom(u_f)=\Cal V$ fixed) 
are continuous (with respect to the topology defined in \S 4.0).  
\endproclaim 

The proof is left to the reader.  

For a function $f \in \N'$ with $\alpha(f) \neq 0$ (i.e. $\sigma(f) \neq 0$), 
let us introduce a new coordinate $w \in \C$ by 
$$ z = \t_f (w) \equiv  \frac \s {1-e^{-2\pi i \alpha w}},  \tag{A.1.5}$$ 
where $\s=\s (f)$ and $\alpha=\alpha(f)$.  
And define the map $F_f(w)$ by 
$$
F_f(w) = w + \frac 1{2\pi i\alpha} 
\log \left( 1- \frac {\s u(z)}{1+zu(z)} \right) \quad \text{ with } z=\t_f(w) 
\tag{A.1.6}
$$
and 
$$T_f(w)=w-\frac 1{\alpha} .$$
Here $F_f(w)$ is defined for $w$ such that  
$|\s u(z)/(1+zu(z))| \leq 1/2$, 
and the above formula defines a single-valued function using the branch of 
logarithm with 
$-\pi <\Im \log (\cdot) \leq \pi$.

\proclaim{Lemma A.1.7} {\rm (Properties of $F_f$, $\t_f$ and $T_f$)}
There exist $R_0>0$ and a neighbourhood $\N \subset \N'$ of $f_0$ such that 
if $f \in \N$ and $\alpha(f) \neq 0$, then: 
\newline 
{\rm (i)} The map $\t_f:\C \to \CC-\{0, \s(f)\}$, $w \mapsto z=\t_f(w)$ 
is a universal covering, whose covering transformation group is generated by 
$T_f$;  
$\t_f(w) \to 0$ as $\Im \alpha w \to \infty$, and 
$\t_f(w) \to \s$ as $\Im \alpha w \to -\infty$;   
\newline   
{\rm (ii)} If 
$$w \in \C - \bigcup _{n \in \Z} T_f^n D_{R_0}, \ \text{ where } 
D_{R_0}=\{w'\bigm| |w'|<R_0\}, $$ 
then  $\t_f(w) \in \Cal V$ and $|\s u(z)/(1+zu(z))| \leq 1/2$, 
hence (A.1.6) is well-defined and moreover satisfies
$$
|F_f(w)-(w+1)|<\frac 14 , \ |F'_f(w)-1|< \frac 14    
$$  

Moreover $F_f(w) = w+1 + O(1/w^2)$ as $\Im \alpha w \to 0$;   
\smallskip 
\noindent 
{\rm (iii)} $f \circ \t_f = \t_f \circ F_f$  and  $T_f \circ F_f = F_f \circ T_f$;  
\smallskip 
\noindent 
{\rm (iv)} When $f \to f_0$, $\t_f(w) \to \t_0(w)$  uniformly on 
$\{w\bigm| |\Re \alpha|<\frac 34 \}$ and 
$F_f(w)  \to F_0(w)$  uniformly on 
$\C - \cup _{n \in \Z} T_f^n D_{R_0}$.

\endproclaim 
The proof is immediate by a computation and left to the reader.

\bigskip 

\subheading{A.2 General construction}  

For $b_1, b_2 \in \C$ with $\Re b_1 < \Re b_2$, define 
$$ 
\Cal Q(b_1, b_2) =\{z \in \C| \ \Re (z-b_1)> -|\Im (z-b_1)|, \ 
\Re (z-b_2) < |\Im (z-b_2)|\ \}.  
$$
If $b_1=-\infty$ (resp. $b_2=\infty$), the condition involving 
$b_1$ (resp. $b_2$) should be removed.

\proclaim{Proposition A.2.1} 
Let $F$ be a holomorphic function defined in $\Q=\Q(b_1,b_2)$, 
where $\Re b_2 > \Re b_1 +2$  (here $b_1$ or $b_2$ may be $-\infty$ or 
$\infty$).  Suppose 
$$|F(z)-(z+1)|< \frac14, \ \text{ and } |F'(z)-1| < \frac 14 
\ \text{ for } z \in \Q. \tag {A.2.2}$$
Then
\newline
{\rm (0)} $F$ is univalent on $\Q$.  
\newline
{\rm (i)} Let $z_0 \in \Q$ be a point such that  $\Re b_1 < \Re z_0 <\Re b_2 -5/4$.  
Denote by $S$ the closed region (a strip) bounded by the two curves 
$\ell=\{z_0+iy| y \in \R\}$ and $F(\ell)$.  Then for any $z \in \Q$, 
there exists a unique $n \in \Z$ such that $F^n(z)$ is defined and 
belongs to $S-F(\ell)$.  
\newline
{\rm (ii)} There exists a univalent function $\Phi : \Q\to \C$
satisfying 
$$\Phi(F(z))=\Phi(z) + 1 $$
whenever both sides are defined.  
Moreover $\Phi$ is unique up to addition of a constant.  
\newline
{\rm (iii)} If we normalize $\Phi$ by $\Phi(z_0)=0$, where $z_0 \in \Q$, 
then the correspondence $F \mapsto \Phi$ is continuous 
with respect to the compact-open topology.  (See also \S 4.0 Notations.)  
 
\endproclaim 

\demo{Proof}  
(0) and (i) are easy and left to the reader.  
\newline 
(ii) Let $z_0 \in \Q$ be as in (i).  
Define $h_1: \{z|\ 0 \leq \Re z \leq 1\} \to \Q $ by 
$$h_1(x+iy)=(1-x)(z_0+iy) + xF(z_0+iy) , \ \text{ for } 0\leq x\leq 1,\ y \in\R.$$
Then we have 
$$
\frac {\partial h_1}{\partial x}=F(z_0+iy)-(z_0+iy), \ 
\frac {\partial h_1}{\partial y}=ixF'(z_0+iy)+i(1-x).  
$$
Hence 
$$
\biggl|\frac {\partial h_1}{\partial z}-1 \biggr|
=\frac12 \bigl|\{F(z_0)-(z_0+iy+1)\}+x(F'(z_0+iy)-1)\bigr| \leq \frac 14, 
$$
$$
\biggl|\frac {\partial h_1}{\partial \bar z}\biggr|
=\frac12 \bigl|\{F(z_0)-(z_0+iy+1)\}-x(F'(z_0+iy)-1) \bigr| \leq \frac 14.  
$$ 
 
Therefore $\bigl| \frac {\partial h_1}{\partial \bar z}/
\frac {\partial h_1}{\partial z}\bigr|<1/3$ and 
$h_1$ is a quasiconformal mapping onto the strip $S$, and satisfies 
$h_1^{-1}(F(z))=h_1^{-1}(z)+1$  for  $ z \in \ell$.   
Let $\sigma_0$ be the standard conformal structure of $\C$, and take the 
pull-back $\sigma = h_1^*\sigma_0$ on $ \{z|\ 0 \leq \Re z \leq 1\}$.  
Then extend $\sigma$
to $\C$ by  $\sigma=(T^n)^* \sigma$ on $\{z| -n \leq \Re z \leq -n+1\}$, 
where $T(z)=z+1$ .  
By Ahlfors-Bers measurable mapping theorem [A], there exists a unique 
quasiconformal mapping $h_2:\C \to \C$ such that $h_2^*\sigma_0=\sigma$ 
and $h_2(0)=0, \ h_2(1)=1$.  By the definition of $\sigma$,  $T$ preserves 
$\sigma$.  Hence $h_2 \circ T \circ h_2^{-1}$ preserves the standard conformal 
structure $\sigma_0$, i.e. conformal, therefore it must be an affine function.  
Since  $T$  has no fixed point in $\C$, so does  $h_2 \circ T \circ h_2^{-1}$, 
hence it is a translation.  Using  $h_2 \circ T \circ h_2^{-1}(0)=1$, we 
have  $h_2 \circ T \circ h_2^{-1}=T$.  

  Now define $\Phi$ by 
$\Phi=h_2 \circ h_1^{-1} $ on $S$ , and extend to the whole $\Q$ using 
the relation 
$ \Phi(F(z))=\Phi(z)+1 .$
Then $\Phi$ is well-defined by (i), continuous and homeomorphic by the 
above relations on $h_1^{-1}$ and $h_2$.   
Moreover $\Phi$ is analytic outside the orbit of $\ell$, 
then analytic in the whole $\Q$ by Morera's theorem.  Thus we have obtained 
the desired univalent function $\Phi$.

If $\Phi'$ is another such function, 
then $\Phi''(z)=\Phi' \circ \Phi^{-1}$ commutes with $T$ at least in $\Phi(\Q)$.  
Hence $\Phi''(z)-z$ extends to $\C$ as a periodic function, then $\Phi''$ extends 
to $\C$ as a holomorphic function commuting with $T$.  Similarly ${\Phi''}^{-1}$ 
also has this property, therefore $\phi$ must be an affine function.  However 
an affine function commuting with the translation $T$ is also a translation by 
a constant.  Hence the assertion follows.  

(iii) Let us consider $F$ and $F_0$ defined in the same $\Q$ and satisfying the 
condition of the Proposition.  As in (ii), 
we can construct $h_1$, $h_2$, $\Phi$ for $F$ and $h_{1,0}$, $h_{2,0}$, 
$\Phi_0$ for $F_0$.  
It is easy to see that on any compact set of $\{z| 0\leq \Re z \leq 1\}$, 
$\frac {\partial h_1}{\partial \bar z}/
\frac {\partial h_1}{\partial z} \to 
\frac {\partial h_{1,0}}{\partial \bar z}/
\frac {\partial h_{1,0}}{\partial z}$ 
as $F \to F_0$.   Hence 
$\sigma_F=h_1^* \sigma_0 \to \sigma_{F_0}=h_{1,0}^* \sigma_0$ and 
$h_2 \to h_{2,0}$ on any compact set as $f \to F_0$.  
It follows from the definition of the extension of $\Phi$ that 
$\Phi \to \Phi_0$ as $F \to F_0$.  
\qed \enddemo

  A strip $S$ as in Proposition A.2.1 (ii), is called a {\it fundamental region} 
for the map $F|_\Q$.  The quotient space 
$$\aligned 
\Cal C &=S/_\sim, \ \text{ where } \ell \ni z \sim F(z) \in F(\ell) \\ 
&= \Q/_\sim, \ \text{ where }  z \sim F(z) \text{ if } 
z \in \Q \cap F^{-1}(\Q)  
\endaligned $$
is topologically a cylinder which is called the {\it Ecalle cylinder}.  
Moreover, 
$\Cal C$ has a natural structure of a Riemann surface, using $F$ near $\ell$ 
as a coordinate patching.

\proclaim{Lemma A.2.3} Let $F$, $\Q$, $S$ be as above.  Then 
$\pi \circ \Phi$ induces an isomorphism 
$$\overline \Phi:\Cal C=S/_\sim  \to \C^*=\C-\{0\}.$$
\endproclaim
Proof.  It is easy to see from the construction that $\overline \Phi$ 
is a covering map and induces an isomorphism between the fundamental groups. 
\qed 

\proclaim{Lemma A.2.4}
Suppose that $\Phi$ and $v$ are holomorphic 
functions in a region $\Cal U$ satisfying: 
$$\Phi \ \text{ is univalent in } \Cal U, \ \ 
|v(z)-1|<1/4 \ \text{ for } z \in \Cal U \ \text{ and }$$
$$\Phi(z+v(z))=\Phi(z)+1 , \ \text{  if }\ z, z+v(z) \in \Cal U.$$
\newline 
{\rm (i)}  There exist universal constants $R_1, C_1, C_2>0$  such that if 
$\Cal U=\{z \bigl| \ |z-z_0|<R\}$ for $R \geq R_1$, then 
$$ \bigl|\Phi'(z_0)-\frac 1{v(z_0)}\bigr| \leq 
C_1\left(\frac1{R^2}+|v'(z_0)|\right)
\leq \frac{C_2}{R} .  $$
\newline
{\rm (ii)}  Suppose $\Cal U=\{z \in \C^*| \ \theta_1< \arg z < \theta_2 \}$ 
$(\theta_2 < \theta_1+2 \pi)$ and $|v'(z)|\leq C/|z|^{1+\nu}$ $(z \in \Cal U)$ 
for some $C, \nu >0$.  For $z_0 \in \Cal U$ and $\theta'_1, \ \theta'_2$ 
with $\theta_1<\theta'_1<\theta'_2<\theta_2$, there exists $R_2, C_3>0$ and 
$\xi \in \C$ such that 
$$\bigl|\Phi(z)-\int_{z_0}^z \frac {d\zeta}{v(\zeta)}-\xi \bigr| \leq 
C_3\left(\frac 1{|z|}+\frac C{|z|^\nu}\right),$$
for $z$ satisfying $\theta'_1 < \arg z < \theta'_2$, $dist(z, \C-\Cal U)>R_1$.  
Moreover $C_3$ depend only on $\theta_i, \theta_i'$.  
\endproclaim 

See [Y] for a similar estimate.  

\demo{Proof} (i) We may suppose $z_0=0$.  
We take $R$ so that $R>> 1$.  
It follows from Koebe's distortion theorem [P] that if 
$|z|<R-2$, then 
$$
\frac{|v(z)|}{(1+|v(z)|/2)^2} \leq 
\biggl| \frac{\Phi(z+v(z))-\Phi(z)}{\Phi'(z)} \biggr| \leq 
\frac{|v(z)|}{(1-|v(z)|/2)^2},
$$
since $\Phi$ is univalent in $\{\zeta \bigl| \ |\zeta -z|<2\ \}$.   Hence 
$C \leq |\Phi'(z)| \leq C'$  if  $|z|<R-2$.  
We have 
$|\Phi''(z)|\leq C/R$  if $|z|<R/2$. 
(In fact, by Cauchy's formula, $\Phi''(z)$ can be expressed in terms of 
an integral of $\Phi'(\zeta)/(\zeta-z)^2$ over the contour 
$\{\zeta \bigl| \ |\zeta-z|=R/2\}$, then use the above estimate.)
Using the formula 
$$
\Phi(z+a)=\Phi(z)+a\Phi'(z)+a^2\int_0^1 (1-t)\Phi''(z+at)dt, 
$$
we obtain $|1-\Phi'(z)v(z)| \leq C/R$ if  $|z|<R/2-5/4$.  
Again by Cauchy's formula, 
$|(1-\Phi'(z)v(z))'|=|\Phi''(z)v(z)+\Phi'(z)v'(z)| \leq C/R^2$  
if $|z|<R/4$.  It also follows from Cauchy's formula that 
$|v''(z)|\leq C/R^2$ for $|z|<R/2$.  Hence we have 
$|\Phi''(z)| \leq C(1/R^2+|v'(0)|)$ for $|z|<R/4$.  By the above formula 
again, we have $|1-\Phi'(0)v(0)|\leq C(1/R^2+|v'(0)|)$.

\bigskip
\enddemo
\subheading{A.3.  Fatou coordinates for $f_0$---Proof of 4.1}

Let $\Q_0^+=\Q(\xi_1, \infty)$, $\Q_0^-=\Q(-\infty, -\xi_1)$ for $\xi_1>0$.  
If $\xi_1$ is large enough, $F_0$ satisfies (A.2.2) on $\Q_0^+$ and $\Q_0^-$.  
Hence by Proposition A.2.1, there exist univalent analytic functions 
$\Phi_{+,0}:\Q_0^+ \to \C$ and $\Phi_{-,0}:\Q_0^- \to \C$ satisfying 
$\Phi_{+,0}(F_0(w))=\Phi_{+,0}(w)+1$ and $\Phi_{-,0}(F_0(w))=\Phi_{-,0}(w)+1$.  
If 
$$F_0(w)=w+1+\frac aw +O(\frac 1{w^2}),$$ 
then by Lemma A.2.4 (ii), we have 
$$
\Phi_{\pm,0}(w)=w - a \log w + c_{\pm,0} +o(1)
$$ 
as $w$ tend to $\infty$ within a sector as in Lemma A.2.4 (ii) with 
$\Cal U=\Q_0^\pm$, 
where $c_{\pm,0}$ are constants and 
the branches of logarithm for $\Phi_{\pm,0}$ are chosen so that 
they coincide in the upper component of $\Q_0^+ \cap \Q_0^-$ and differ 
by $2 \pi i$ in the lower.  
Then for large $\xi_0>0$, $\Q_0$ as in (4.1.2) is contained in $\Phi_{-,0}(\Q_0^-)$
Define 
$\Omega_+=\t_0(T(\Q_0^+))$, $\Omega_-=\t_0(\Phi_{-,0}^{-1}(\Q_0))$.  
The properties in (4.1.1) are easily verified.  (If necessary, take a larger 
$\xi_0$.)  

Let $\phi_0=\t_0 \circ \Phi_{-,0}^{-1}$ hence 
$\phi^*_0=\t_0^{-1} \circ \phi_0=\Phi_{-,0}^{-1}$.  Then by the above, 
$\Q_0 \subset \Dom(\phi_0)$ and $\phi_0(w+1)=f_0 \circ \phi_0(w)$  
if $w, w+1 \in \Q_0$.  Using this relation, we extend $\phi_0$ 
to the maximal domain 
$\{ w \in \C \bigm| \ w-n \in \Q_0 \text{ for an integer } n \geq 0 
\text{ and } f_0^j(\phi_0(w-n)) \in \Dom(f_0) \text{ for } j=0, \dots, n-1 \}$.  
Note that for a large $\eta_0>0$, $\{w \bigm| |\Re w| \leq 1/2,\ 
|\Im w|>\eta_0\}$ is contained in $\Q_0$ and 
its image by $\phi_0$ is contained in $\t_0(\Q_0^+) \subset \Cal B$.  
The rest of (4.1.2) can be checked easily.  

Let $\Phi_0=\Phi_{+,0} \circ \t_0^{-1}$ and we extend it similarly to $\Cal B$ 
using the functional equation for $\Phi_{+,0}$.  Hence (4.1.3) follows.  

The properties of $\tilde {\Cal B}$ in (4.1.4) follow from the above.  
The function $\tE_{f_0}=\Phi_0 \circ \phi_0=\Phi_{+,0} \circ \Phi_{-,0}^{-1}$ 
satisfies 
the functional equation and is univalent in $\{w \bigm| |\Im w|> \eta'_0\}$ for 
a large $\eta'_0>0$.  Moreover $\Im \tE_{f_0}(w) \to \pm \infty$ as 
$\Im w \to \pm \infty$.  Hence $\cE_{f_0}$ can be extended to $0$ and $\infty$ 
conformally.  This proves (4.1.4).  

By the normalization (4.1.5), we have $c_{+,0}=c_{-,0}$.

\bigskip
\subheading{A.4.  Fatou coordinates for $f \in \N \cap \Cal F_1$---
Proof of 4.2} 

Let $f \in \N \cap \Cal F_1$, where $\Cal F_1$ is as in \S 4.2.  
Define $\Q_f^+=\Q(\xi_1, -\xi_1+\frac 1\alpha)$ and 
$\Q_f^-=\Q(\xi_1-\frac 1\alpha, -\xi_1)=T_f(\Q_f^+)$ for $\xi_1>0$.  
One can choose $\N$ small and $xi_1$ large so that 
$\Re \frac 1\alpha>2 \xi_1 + 2$, $\Q_f^\pm \subset \Dom(F_f)$ and 
$F_f$ satisfies (A.2.2) on 
$\Q_f^\pm$.  By Proposition A.2.1, there exist univalent analytic functions 
$\Phi_{+,f}:\Q_f^+ \to \C$ and $\Phi_{-,f}:\Q_f^- \to \C$ satisfying 
$\Phi_{+,f}(F_f(w))=\Phi_{+,f}(w)+1$ and $\Phi_{-,f}(F_f(w))=\Phi_{-,f}(w)+1$.  
We fix two points $w_+ \in \Q_0^+$ and $w_- \in \Q_0^-$, and  normalize 
$\Phi_{\pm,f}$ by setting $\Phi_{\pm,f}(w_\pm)=\Phi_{0,f}(w_\pm)$.

Since $F_f(w)=w+1 + O(1/w^2)$ as $\Im \alpha w \to 0$ by Lemma A.1.7, 
it follows from Lemma A.2.4 (ii) that there exist constants $c_\pm=c_\pm(f)$ 
such that 
$$\Phi_{\pm,f}(w)=z+c_\pm +o(1) \tag{A.4.1}$$ 
when $w$ tend to $\infty$ within a sector 
of the form 
$\{ w \bigm| \theta'_1< \arg(w-\omega_0) <\theta'_2 \} \subset \Q_f^\pm$, 
where $1/3<\theta'_1<theta'_2<3/4$ and $\omega_0 \in \C$.

The function $\Phi_{-,f} \circ T_f$ is defined on $\Q_f^+$ 
and satisfies the same functional equation as $\Phi_{+,f}$ by Lemma A.1.7 (iii).  
It follows from the uniqueness in Proposition A.2.1 (ii)  that there exists
a constant $d_1=d_1(f)$ such that 
$$\Phi_{-,f} \circ T_f = \Phi_{+,f} + d_1. \tag{A.4.2}$$  
In fact, (A.4.1) determines the constant: $d_1=c_-(f) -c_+(f) -\frac 1\alpha$.  

Now define $\tE_f(w)=\Phi_{+,f} \circ \Phi_{-,f}^{-1}$ on 
$\Phi_{-,f}(\Q_f^+ \cap \Q_f^-)$, which contains at least the vertical strip 
$\{w \bigm| |\Re w| \leq 1,\ |\Im w| \geq \eta_0\}$ for a large $\eta_0>0$.  
This $\eta_0$ can be chosen uniformly for $f\in \N \cap \Cal F_1$ near $f_0$.  
It satisfies the functional equation $\tE_f(w+1)=\tE_f(w)+1$ whenever the 
both sides are defined.  By this relation, $\tE_f$ can be 
extended to $\{w \bigm| |\Im w|>\eta_0\}$.  By (A.4.1), we have 
$$\tE_f(w)=w+c_+(f)-c_-(f)+o(1) \ \text{ when } \Im w \to \infty.  \tag{A.4.3}$$    
Similarly, $\tE_f(w)-w$ tends to a constant as $\Im w \to -\infty$.  

By Proposition A.2.1, we have $\Phi_{+,f} \to \Phi_{+,0}$ and 
$\Phi_{-,f} \to \Phi_{-,0}$ as $f \to f_0$, therefore 
$\tE_f \to \tE_{f_0}$ uniformly on $\{w\bigm| 
0\leq \Re w \leq 1, \ \Im w = \eta_0\}$, then 
$$c_+(f)-c_-(f)=\int_{i\eta_0}^{1+i\eta_0}(\tE_f(w)-w)dw \to 
\int_{i\eta_0}^{1+i\eta_0}(\tE_{f_0}(w)-w)dw=
c_{+,0}-c_{-,0}=0 , 
\tag{A.4.4}$$ 
where the integrals are over the segment joining $i\eta_0$ and $1+i\eta_0$.

\medskip 
Let us show that $\Phi_{-,0}(\Q_f^-)$ contains $\Q_f$ as in \S 4.2, for 
large $\xi_0, \eta_0>0$.  Combining (A.1.6), Lemma A.2.4 and Proposition A.2.1 (iii), 
one can show that the right boundary curve of $\Q_f$ is contained in 
$\Phi_{-,0}(\Q_f^-)$, if $\xi_0$ and $\eta_0$ are large enough.  Similarly, 
The left boundary curve of $T_f^{-1}(\Q_f)$ is contained in $\Phi_{+,0}(\Q_f^+)$.  
Using $\Q_f^-=T_f(\Q_f^+)$, (A.4.2) and (A.4.4) and increasing $\xi_0, \eta_0$ 
if necessary, we conclude that the left boundary curve of $\Q_f$ is contained in 
$\Phi_{-,0}(\Q_f^-)$.  Therefore $\Q_f \subset \Phi_{-,0}(\Q_f^-)$.  

\medskip 

Define $\phi_f=\t_f \circ \Phi_{-,f}^{-1}$, then we have  
$\Q_f \subset \Dom(\phi_f)$ and $\phi_f(w+1)=f \circ \phi_(w)$  
if $w, w+1 \in \Q_f$.  Using this relation, extend $\phi_f$ 
to the domain 
$\{ w \in \C \bigm| \ w-n \in \Q_f \text{ for an integer } n \geq 0 
\text{ and } f^j(\phi_f(w-n)) \in \Dom(f) \text{ for } j=0, \dots, n-1 \}$.  
Then $\phi_0$ satisfies (4.2.1).  

Let $\tR_f=\Phi_{-,f} \circ T_f \circ \phi_f$. Then $\tR_f$ coincides 
with $\tE_f+c_-(f)-c_+(f)-\frac 1\alpha$ on its domain of definition and 
satisfies $\tR_f(w+1)=\tR_f(w)+1$.  So $cR_f=\pi \circ \tR_f \circ \pi^{-1}$ 
is well-defined and extends analytically to $0$ by (A.4.3) and 
similarly to $\infty$.  
Since $\tR_f(w)=w-\frac 1\alpha +o(1)$ as $\Im w \to \infty$, 
$\cR'_f(0)=\exp(-2\pi i \frac 1\alpha)$.  Thus (4.2.2) is proved.  

\medskip 

It is easy to see that $\phi_f \to \phi_0$ and $\tE_f \to \tE_{f_0}$ 
when $f \in \N \cap \Cal F_1$ and $f \to f_0$.  Hence 
$$\tR_f + \frac 1 \alpha = \tE_f +(c_-(f)-c_+(f)) \to \tE_{f_0}.  $$  
Using the fact that $cR_f$ is extended analytically to $0$ and $\infty$,  
we have 
$$e^{2\pi i/\alpha}\cR_f \to \cE_{f_0}.$$
Thus (4.2.4) is proved.  

In order to see (4.2.3), let us first verify (4.3.4).  In fact the latter is 
immediate from the definition of $\phi_f$ and $\tR_f$, $\t_f \circ T_f = \t_f$ 
and the functional equation for $\phi_f$.  
To obtain (4.2.3), we restrict $\Dom(\tR_f)$ so that 
$$|\tR_f(w)+\frac 1\alpha-w|<\frac 1{2\alpha} \text{ on } \Dom(\tR_f). 
\tag{A.4.5}$$  
Note that (4.2.4) is still true.  
Now suppose that 
$w, w' \in \Dom(\phi_f)$ and $\cR_f^m(\pi(w))=\pi(w')$ for a positive integer 
$m$.  Then for the lift, there exists an integer $n$ such that 
$\tR_f^m(w)+n=w'$.  By (A.4.5), we have $n>m$.  
Hence the assertions of (4.2.3) follow from (4.3.4).  

The corollary (4.2.5) follow from (4.1.2) and (4.2.4).  

\medskip 
One can prove (4.3.1) similarly.  

\bigskip 

\subheading{A.5.  Parabolic fixed points with the multiplier $\neq 1$}  

We only state the coordinate changes which correspond to A.1.  
The rest of arguments in A.2-A.4 is immediately generalized to this case.  

\medskip 

First, let $f_0$ be as in \S 7.  The exists a coordinate near $0$, in which 
$f_0^q$ has the form 
$$f_0^q(z)=z+z^{q+1}+O(z^{2q+1}). $$
Let us introduce a coordinate $w$ by 
$$z=-\frac 1{qw^q}.$$  
Then the corresponding map in this coordinate has the form 
$$F_0(w)=w+1+O(1/w), $$  
where note that $F_0$ is multi-valued.  

The coordinate $w$ should be understood in terms of the Riemann surface 
$\Cal W$ as follows.  Let $\Q_0^{+,k}$, $\Q_0^{-,k}$ $(k \in \Z/q\Z)$ 
be $q$ copies of $\Q_0^+$, $\Q_0^-$ in A.3.  The intersection 
$\Q_0^+ \cap \Q_0^-$ has two components---the upper and lower triangles.  
Identify the upper triangle of $\Q_0^{+,k}$ with that of $\Q_0^{-,k}$;  
and the lower triangle of $\Q_0^{+,k-1}$ with that of $\Q_0^{-,k}$ $(k \in \Z/q\Z)$.  
The obtained Riemann surface $\Cal W$ is isomorphic to a punctured disc  
and the map $w \mapsto z=-1/qw^q$ gives a conformal map map from 
$\Cal W$ onto a punctured neighbourhood of $0$ in $\C$.  
We consider that the map $F_0$ sends $\Q_0^{+,k}$ into $\Q_0^{+,k+p}$, 
and $\Q_0^{-,k}-\text{(a neighbourhood of the boundary curve)}$ into $\Q_0^{-,k+p}$. 

\medskip 

Now let $f \in \Cal F_1$ and $\sigma^{(k)}(f)$ 
the periodic points of $f$ as in \S 7.   We choose a new coordinate near $0$ 
so that $\sigma^{(k)}(f)=e^{2\pi i k/q} \sigma^{(0)}(f)$, as follows.  
Put $$(z-\sigma^{(0)}(f)) \cdots (z-\sigma^{(k-1)}(f)) 
=z^q-b_{q-1}(f)z^{q-1}-\cdots -b_1(f)z-b_0(f).$$    
Then it is easily seen that 
$$b_0(f)=-2\pi i \alpha(f)(1+o(1)) \text{ and } b_j(f)=O(\alpha(f)) 
\ (j=1, \dots, q-1)\ \text{ as } f \to f_0.  $$  
Define a new coordinate $z'$ by 
$$z'=\frac z{\left(1+\frac {b_1}{b_0}z+ \cdots + 
\frac {b_{q-1}}{b_0}z^{q-1}\right)^{1/q}}.$$  
This is a well defined coordinate near $0$, and in this coordinate 
we have  $\sigma^{(k)}=e^{2\pi i k/q} \sigma^{(0)}$.  

So we can write 
$$f^q(z)=z+z(z^q-\sigma(f)^q)u_f(z), $$
where $\sigma(f)=\sigma^{(0)}(f)$ and $u_f(z)$ is a non-zero analytic function 
defined near $0$.  Finally, we introduce the coordinate $w$ by 
$$\frac {z^q}{z^q-\sigma^q}=e^{2 \pi i \alpha w}$$ 
and the map $F_f$ by 
$$F_f(w)=w+\frac 1{2\pi i \alpha} \log 
\left(1- \frac {\s^q u_f(z)}{1+zu_f(z)} \right) .
$$ 
Here the coordinate $w$ should be interpreted on a suitable Riemann surface 
which is isomorphic to a neighbourhood of $0$ in $\C$ with $0$ and 
$\sigma^{(k)}(f)$ removed, 
then the $z$ in the definition of $F_f$ can make sense.  

For these maps, one can obtain analogous results as in A.2-A.4.   
The detail is left to the reader.  

\bigskip 
\bigskip
\bigskip

\comment 
----------------------------------------------------------------------------------
\endcomment

\heading
References
\endheading

\noindent 
[A] Ahlfors, L.,  Lectures on quasiconformal mappings, Van Nostrand, 1966. 
\smallskip
\noindent  
[BC] Benedicks, M., and Carleson, L., The dynamics of the H\'enon map, 
Annals of Math., 133 (1991) p. 73-169.  
\smallskip
\noindent  
[Bl] Blanchard, P.,  Complex Analytic dynamics on the Riemann sphere, Bull. 
AMS 11 (1984), p. 85-141.  
\smallskip \noindent  
[BR] Bers, L. and Royden, H.L., Holomorphic families of injections, 
Acta Math., 157 (1986) p. 259-286.  
\smallskip \noindent  
[Bo1] Bowen, R., Equilibrium states and the ergodic theory of 
Anosov diffeomorphisms,  Lecture note in Math. no. 470. Springer, Berlin, 1975.  
\smallskip \noindent  
[Bo2] -----, Hausdorff dimension of quasi-circles, Publ. Math. I.H.E.S. 50 (1979), 
P.11-26.  
\smallskip \noindent  
[DH] Douady, A. and Hubbard, J. H., \'Etude dynamique des polyn\^omes 
complexes, Publ. Math. d'Orsay, 1er partie, 84-02;  2me partie, 85-04.  
\smallskip \noindent  
[J] Jakobson, M.,  Absolutely continuous invariant measures for one-parameter 
families of one-dimensional maps, Comm. Math. Phys., 81 (1981) p. 39-88.  
\smallskip \noindent  
[La] Lavaurs, P., Syst\`emes dynamiques holomorphes: explosion de points 
periodiques paraboliques, Th\`ese de doctrat de l'Universit\'e de Paris-Sud, 
Orsay, France, 1989.  
\smallskip \noindent  
[Ly1] Lyubich, M.Yu., An analysis of stability of the dynamics of rational 
functions, Selecta Mathematica Sovietica, 9 (1990) p.69-90.  
\smallskip \noindent  
[Ly2] -----, On the Lebesgue measure of the Julia set of a quadratic polynomial, 
Preprint SUNY Stony Brook, 1991.  
\smallskip \noindent  
[Ma] Mandelbrot, B.,  On the dynamics of Iterated maps V: Conjecture that 
the boundary of the M-set has a fractal dimension equal to 2,  p. 235-238, 
Choas, Fractals and Dynamics, Eds. Fischer and Smith, Marcel Dekker, 1985.  
\smallskip \noindent  
[Mc] McMullen, C., Area and Hausdorff dimension of Julia sets of entire 
functions, Trans. AMS, 300 (1987) p.329-342.  
\smallskip \noindent  
[Mi] Milnor, J., Dynamics in one complex variables: Introductory lectures, 
Preprint SUNY Stony Brook, Institute for Mathematical 
Sciences, 1990.  
\smallskip \noindent  
[Mi2] -----, Self-similarity and hairiness in the Mandelbrot set, 
in: Computers in Geometry and Topology, ed. Tangora, M.C, Lect. Notes 
in Pure and Appl. Math., Vol. 114, 1989, Marcel Dekker, p.211-257.  
\smallskip \noindent  
[MSS] Ma\~n\'e, R., Sad, P.,  and Sullivan, D., On the dynamics of 
rational maps, Ann. scient. Ec. Norm. Sup., (4) 16 (1983) 
p.193-217. 
\smallskip \noindent 
[P] Pommerenke, Ch., Univalent functions, Vandenhoeck \& Ruprecht, G\"ottingen, 
1975. 
\smallskip \noindent    
[Re] Rees, M., Positive measure sets of ergodic rational maps, 
Ann. scient. Ec. Norm. Sup., (4) 19 (1986) p. 383-407.  
\smallskip \noindent  
[Ru] Ruelle, D.,  Repellers for real analytic maps, Erg. Th. \& Dynam. Sys. 
2 (1982), p.99-107.  
\smallskip \noindent 
[Sh] Shishikura, M., in preparation.   
\smallskip \noindent  
[ST] Sullivan, D. and Thurston, W.P., Extending holomorphic motions, 
Acta Math., 157 (1986) p. 243-257.  
\smallskip \noindent  
[Su] Sullivan, D., Conformal dynamical systems, in: Geometric dynamics, 
Lecture note in Math. 1007 (1983) p.725-752, Springer, Berlin.  
\smallskip \noindent  
[T] Tan Lei, Similarity between the Mandelbrot set and Julia sets, 
Commun. Math. Phys., 134 (1990) p.587-617.  
\smallskip \noindent  
[Y] Yoccoz, J.-Ch., Lin\'earisation des germes de diff\'eomorphismes holomorphes 
de $(\C, 0)$, C. R. Acad. Sci. Paris, t.306, S\'er. I (1988) p.55-58.   

\comment

\newline 
[EL] Eremenko, A.E. and Lyubich, M.Yu.,  The dynamics of analytic 
transformations, Leningrad Math.J, 1 (1990) p. 563-634.   

[H] Herman, M. R., 

[DD]? Devaney, R. and Douady, A., 
\newline 

[HKS] Handler, I., Kauffman, L.H. and Sandin, D., On Crossing the boundary 
of the Mandelbrot set, see [Mi2], p.151-177.  
\endcomment

\bigskip 
\bigskip 

\noindent 
Address: 

\medskip 
\noindent 
Tokyo Institute of Technology 
\newline  
Department of Mathematics 
\newline 
Ohokayama, Meguro, Tokyo 152, Japan 
\newline  
email: mitsu\@math.titech.ac.jp 

\smallskip 
\noindent 
and 

\smallskip 
\noindent 
Institute for Mathematical Sciences 
\newline 
State University of New York at Stony Brook 
\newline 
Stony Brook, NY 11794-3660, USA 
\newline 
email: mitsu\@math.sunysb.edu

\end